\documentstyle{article}
\catcode`\@=11
\def\marginnote#1{}

\newcount\hour
\newcount\minute
\newtoks\amorpm
\hour=\time\divide\hour by60
\minute=\time{\multiply\hour by60 \global\advance\minute by-\hour}
\edef\standardtime{{\ifnum\hour<12 \global\amorpm={am}%
        \else\global\amorpm={pm}\advance\hour by-12 \fi
        \ifnum\hour=0 \hour=12 \fi
        \number\hour:\ifnum\minute<10 0\fi\number\minute\the\amorpm}}
\edef\militarytime{\number\hour:\ifnum\minute<10 0\fi\number\minute}

\def\draftlabel#1{{\@bsphack\if@filesw {\let\thepage\relax
   \xdef\@gtempa{\write\@auxout{\string
      \newlabel{#1}{{\@currentlabel}{\thepage}}}}}\@gtempa
   \if@nobreak \ifvmode\nobreak\fi\fi\fi\@esphack}
        \gdef\@eqnlabel{#1}}
\def\@eqnlabel{}
\def\@vacuum{}
\def\draftmarginnote#1{\marginpar{\raggedright\scriptsize\tt#1}}

\def\draft{\oddsidemargin -.5truein
        \def\@oddfoot{\sl preliminary draft \hfil
        \rm\thepage\hfil\sl\today\quad\militarytime}
        \let\@evenfoot\@oddfoot \overfullrule 3pt
        \let\label=\draftlabel
        \let\marginnote=\draftmarginnote
   \def\@eqnnum{(\theequation)\rlap{\kern\marginparsep\tt\@eqnlabel}%
\global\let\@eqnlabel\@vacuum}  }
\newcounter{app}
\newcounter{sapp}[app]

\newcommand{\app}[1]{
\refstepcounter{app}{\noindent\Large\bf Appendixes
 \ #1 \par \vspace{5mm}}
\setcounter{equation}{0}
\def\theequation{\Alph{app}.\arabic{equation}}}
\def\thesapp{\Alph{app}.\arabic{sapp}}
\newcommand{\sapp}[1]{\par \refstepcounter{sapp}{\noindent\large\bf \thesapp
\ #1 \par \vspace{3mm}}
\def\theequation{\Alph{app}.\arabic{equation}}}

\makeatletter
\newdimen\normalarrayskip              % skip between lines
\newdimen\minarrayskip                 % minimal skip between lines
\normalarrayskip\baselineskip
\minarrayskip\jot
\newif\ifold             \oldtrue            
\def\arraymode{\ifold\relax\else\displaystyle\fi} % mode of array entries
\def\eqnumphantom{\phantom{\mbox{\rm
(\theequation)}}}% right phantom in eqnarray
\def\@arrayskip{\ifold\baselineskip\z@\lineskip\z@
     \else
     \baselineskip\minarrayskip\lineskip2\minarrayskip\fi}
\def\@arrayclassz{\ifcase \@lastchclass \@acolampacol \or
\@ampacol \or \or \or \@addamp \or
   \@acolampacol \or \@firstampfalse \@acol \fi
\edef\@preamble{\@preamble
  \ifcase \@chnum
     \hfil$\relax\arraymode\@sharp$\hfil
     \or $\relax\arraymode\@sharp$\hfil
     \or \hfil$\relax\arraymode\@sharp$\fi}}
\def\@array[#1]#2{\setbox\@arstrutbox=\hbox{\vrule
     height\arraystretch \ht\strutbox
     depth\arraystretch \dp\strutbox
     width\z@}\@mkpream{#2}\edef\@preamble{\halign \noexpand\@halignto
\bgroup \tabskip\z@ \@arstrut \@preamble \tabskip\z@ \cr}%
\let\@startpbox\@@startpbox \let\@endpbox\@@endpbox
  \if #1t\vtop \else \if#1b\vbox \else \vcenter \fi\fi
  \bgroup \let\par\relax
  \let\@sharp##\let\protect\relax
  \@arrayskip\@preamble}
\def\eqnarray{\stepcounter{equation}%
              \let\@currentlabel=\theequation
              \global\@eqnswtrue
              \global\@eqcnt\z@
              \tabskip\@centering
              \let\\=\@eqncr
              $$%
 \halign to \displaywidth\bgroup
    \eqnumphantom\@eqnsel\hskip\@centering
    $\displaystyle \tabskip\z@ {##}$%
    &\global\@eqcnt\@ne \hskip 1.2\arraycolsep
         $\displaystyle\arraymode{##}$\hfil
    &\global\@eqcnt\tw@ \hskip 1.2\arraycolsep
         $\displaystyle\tabskip\z@{##}$\hfil
         \tabskip\@centering
    &{##}\tabskip\z@\cr}
\makeatother
\catcode`\@=11
\ifcase\@ptsize
 \font\tenmsa=msam10
 \font\sevenmsa=msam7
 \font\fivemsa=msam5
 \font\tenmsb=msbm10
 \font\sevenmsb=msbm7
 \font\fivemsb=msbm5
 \font\teneu=eufm10
 \font\seveneu=eufm7
 \font\fiveeu=eufm5
 \font\tenib=cmmib10
 \font\sevenib=cmmib7
 \font\fiveib=cmmib5
\or
 \font\tenmsa=msam10 scaled \magstephalf
 \font\sevenmsa=msam7 scaled \magstephalf
 \font\fivemsa=msam5 scaled \magstephalf
 \font\tenmsb=msbm10 scaled \magstephalf
 \font\sevenmsb=msbm7 scaled \magstephalf
 \font\fivemsb=msbm5  scaled \magstephalf
 \font\teneu=eufm10  scaled \magstephalf
 \font\seveneu=eufm7  scaled \magstephalf
 \font\fiveeu=eufm5   scaled \magstephalf
 \font\tenib=cmmib10  scaled \magstephalf
 \font\sevenib=cmmib7  scaled \magstephalf
 \font\fiveib=cmmib5   scaled \magstephalf
\or
 \font\tenmsa=msam10 scaled \magstep1
 \font\sevenmsa=msam7 scaled \magstep1
 \font\fivemsa=msam5  scaled \magstep1
 \font\tenmsb=msbm10 scaled \magstep1
 \font\sevenmsb=msbm7 scaled \magstep1
 \font\fivemsb=msbm5  scaled \magstep1
 \font\teneu=eufm10   scaled \magstep1
 \font\seveneu=eufm7 scaled \magstep1
 \font\fiveeu=eufm5 scaled \magstep1
 \font\tenib=cmmib10     scaled \magstep1
 \font\sevenib=cmmib7   scaled \magstep1
 \font\fiveib=cmmib5   scaled \magstep1
\fi
\newfam\msafam
\textfont\msafam=\tenmsa  \scriptfont\msafam=\sevenmsa
  \scriptscriptfont\msafam=\fivemsa
\newfam\msbfam
\textfont\msbfam=\tenmsb  \scriptfont\msbfam=\sevenmsb
  \scriptscriptfont\msbfam=\fivemsb
\def\Bbb{\ifmmode\let\next\Bbb@\else
 \def\next{\errmessage{Use \string\Bbb\space only in math mode}}\fi\next}
\def\Bbb@#1{{\Bbb@@{#1}}}
\def\Bbb@@#1{\fam\msbfam#1}
\newfam\eufam
\textfont\eufam=\teneu  \scriptfont\eufam=\seveneu
  \scriptscriptfont\eufam=\fiveeu
\def\frak{\ifmmode\let\next\frak@\else
 \def\next{\errmessage{Use \string\frak\space only in math mode}}\fi\next}
\def\frak@#1{{\frak@@{#1}}}
\def\frak@@#1{\fam\eufam#1}
\newfam\ibfam
\textfont\ibfam=\tenib  \scriptfont\ibfam=\sevenib
  \scriptscriptfont\ibfam=\fiveib
\def\bold{\ifmmode\let\next\bold@\else
 \def\next{\errmessage{Use \string\bold\space only in math mode}}\fi\next}
\def\bold@#1{{\bold@@{#1}}}
\def\bold@@#1{\fam\ibfam#1}
\def\hexnumber@#1{\ifcase#1 0\or 1\or 2\or 3\or 4\or 5\or 6\or 7\or 8\or
 9\or A\or B\or C\or D\or E\or F\fi}
\def\newsymbolb#1#2#3#4{\mathchardef#1="#2\hexnumber@\msbfam#3#4}
\relax
\begin{document}
\def\bea{\begin{eqnarray}}
\def\eea{\end{eqnarray}}
\def\beq{\begin{equation}}          \def\bn{\beq}
\def\eeq{\end{equation}}            \def\ed{\eeq}
\def\nn{\nonumber}                  \def\g{\gamma}
\def\Uq{U_q(\widehat{\frak{sl}}_2)}
\def\Uqp{U_q(\widehat{\frak{sl}}'_2)}
\def\Uqd{U^{*}_q(\widehat{\frak{sl}}_2)}
\def\uq{U_q({sl}_2)}
\def\uqd{U^*_q({sl}_2)}
\def\slaff{\frak{sl}^\prime_2}
\def\aff{\widehat{\frak{sl}}_2}
\def\ot{\otimes}
\def\id{\mbox{\rm id}}
\def\tr{\mbox{\rm tr}}
\def\ctg{\mbox{\rm ctg}}
\def\Re{{\rm Re}\,}
\def\RR{\Bbb{R}}
\def\ZZ{\Bbb{Z}}
\def\CC{\Bbb{C}}
\def\RR{\Bbb{R}}
\def\r#1{\mbox{(}\ref{#1}\mbox{)}}
\def\d{\delta}
\def\D{\Delta}
\def\da{{\partial_\alpha}}
\let\da=p
\def\R{{\cal R}}
\def\h{\hbar}
\def\Ga#1{\Gamma\left(#1\right)}
\def\ep{\varepsilon}
\def\ve{\ep}
\def\fract#1#2{{\mbox{\footnotesize $#1$}\over\mbox{\footnotesize $#2$}}}
\def\stackreb#1#2{\ \mathrel{\mathop{#1}\limits_{#2}}}
\def\res#1{\stackreb{\mbox{\rm res}}{#1}}
\def\lim#1{\stackreb{\mbox{\rm lim}}{#1}}
\def\Res#1{\stackreb{\mbox{\rm Res}}{#1}}
\let\dis=\displaystyle
\def\ee{{\rm e}}
\def\D{\Delta}
\renewcommand{\theequation}{{\thesection}.{\arabic{equation}}}
\def\Y-{\widehat{Y}^-}
\font\fraksect=eufm10 scaled 1728
\def\DYsect{\widehat{DY(\hbox{\fraksect sl}_2)}}
\def\DY{\widehat{DY(\frak{sl}_2)}}
\def\Yd{\DY}
\def\Ydd{\DY}
\def\aa{\alpha}
\def\s{\sigma}
\def\l{\lambda}
\def\la{\land}
\def\L{\Lambda}
\def\ep{\epsilon}
\def\g{\gamma}
\def\d{\delta}
\def\p{\partial}
\def\J{{\bf J(\s)}}
\def\I{{\bf I(\s)}}
\let\z=\alpha
\let\b=\beta
\let\u=u
\let\v=v
\let\w=z
\let\hsp=\qquad

\begin{titlepage}
\begin{center}
\hfill ITEP-TH-8/97\\
\hfill {\tt q-alg/?????}\\
\bigskip\bigskip
{\Large\bf Traces of creation-annihilation operators and Fredholm's formulas}\\
\bigskip
\bigskip
{\large A. Chervov \footnote{E-mail: chervov@vitep1.itep.ru
~~~~ alex@lchervova.home.bio.msu.ru}}\\
\medskip
{\it Faculty of Mathematics and Mechanics  \\
Moscow State University\\}
\medskip{\it and }\\
\medskip
{\it Institute of Theoretical \& Experimental Physics\\
117259 Moscow, Russia}\\
\bigskip
\bigskip
\bigskip
{Revised \today}
\end{center}
\begin{abstract}

We prove the formula for the  traces of certain class of operators
in bosonic and fermionic Fock spaces.
Vertex operators belong to this class. Traces of vertex operators
can be used for calculation of correlation functions and formfactors
of integrable models (XXZ, Sine-Gordon, etc.), that is why we
are interested in this problem. Also we show that Fredholm's minor and
determinant can be expressed by such traces. We obtain a short proof
of the Fredholm's formula for the solution of an integral equation.

\end{abstract}
\end{titlepage}
\clearpage
\newpage

\setcounter{equation}{0}
\setcounter{footnote}{0}

\section{Introduction.}

In this paper  we prove formulas for the traces of some
class of operators in bosonic and fermionic Fock spaces.
Such  traces are used  in mathematical physics,
because  correlators and formfactors
of integrable models can be expressed in  terms of such traces
 (\cite{JM},\cite{L} and \cite{KLP2}). On the other hand
Fredholm's determinant and minor can be expressed
in the terms of similar traces;
and special case of our  formula turns out to be  Fredholm's formula
for the solution of an integral equation.

 Let us discuss the  formulas.
Let $\CC [x_1, x_2,...] $
 be the algebra of polinomials. Let  ${\bf C} (x_1, x_2,...) $ be
the   operator of multiplication on a polinomial $ C (x_1, x_2,...) $;~~ let
${\bf A}(\partial_{x_1},\partial_{x_2},...) $ be a  polinomial of operators
$\partial_{x_1}=\frac{\partial}{\partial{x_1}},$
$\partial_{x_2}=\frac{\partial}{\partial{x_2}},... $;~~ let $\bar\rho$
be an arbitrary, degree preserving
 homomorphism of the algebra  $\CC [x_1, x_2,...] $  onto itself.
We prove the following formula for the trace over the space
$\CC [x_1, x_2,...] $ of the product  of these three
operators:

\bea Tr\left (\bar\rho {\bf A} (\partial_{x_1},\partial_{x_2},....) {
\bf C} (x_1, x_2,...)\right)=
\left (Tr\bar\rho\right)\left < A (\partial_{x_1},\partial_{x_2},....)|
\frac{1}{1 -\bar\rho}C (x_1, x_2,...)\right >\label{tr_1}
\eea
%where
Here pairing $ <...|... > $ is defined by
the formula: $ < (\partial_{x_1}) ^{i_1}... (\partial_{x_n})
^{i_n}|(x_1) ^{j_1}... (x_n) ^{j_n} >=\delta_{i_1}^{j_1}...
\delta_{i_n}^{j_n} i_1!... i_n! $.
The trace of the operator $O$: $\CC [x_1, x_2,...] \to$
$\CC [x_1, x_2,...]  $ is  the sum of diagonal elements
in the natural basis $ (x_1) ^{i_1} (x_2) ^{i_2}....(x_k) ^{i_k}$.

We shall also prove  analogous formula for the anticommuting variables
$\xi_i\xi_j= -\xi_j\xi_i$:
\bea Tr_{\L [\xi_1,\xi_2,...]}\left (\bar\rho {\bf A} (\partial_{\xi_1},
\partial_{\xi_2},....) {\bf C } (\xi_1,\xi_2,...)\right)=
\left (Tr_{\L [\xi_1,\xi_2,...]}\bar\rho\right)\left < A (\partial_{
\xi_1},\partial_{\xi_2},....)|\frac{1}{1+\bar\rho}C (\xi_1,\xi_2,...)
\right>\label{tr_2}\label{first_fml}\eea

The special cases of the first formula can be found in \cite{JM},
special cases of the second one in \cite{SJM}.

These formulas express  matrix elements
of the operator inverse to  $ (1\pm\bar\rho) $. If $\bar\rho$
is an integral operator then the second formula is equivalent
to the  Fredholm's formula for the solution of integral equation.
Thus we obtain  short  proof of the Fredholm's formula, and
 interpretation of Fredholm's minor as the trace over
$\L [\xi_1,\xi_2,...] $ of  some sort of operator, which seems to be
unknown before. (The interpretation
of Fredholm's determinant as such trace is well known.)
Also  our  first formula gives analogous formula of inversion
not via the  Fredholm's determinants, but via Fredholm's permanents
\footnote{permanent of  square matrix is  the  sum of its elements standing in
different columns and different lines}.

Correlation
functions and formfactors of different integrable models can  be found,
as  traces  of vertex operators. This was shown in
\cite{DFJMN}, \cite{JM} for XXZ-model and corresponding six-vertex model and
in \cite{L},\cite{KLP1},\cite{KLP2} and \cite{KLP3} for $SU(2)$-invariant
Thirring model and $Sin-Gordon$ model. The so-called procedure of bosonization
allows one to write down vertex operators, as operators in the
space $\CC [x_1, x_2,...] $ of the  type: $\bar\rho A (
\partial_{x_1},\partial_{x_2},....) C (x_1, x_2,...) $. Our formula
allows to compute the traces of such operators.
Our initial  purpose
was to prove the formula \ref{formul2} and its consequence \ref{calc-gener}
that are used in  \cite{KLP2} to compute
formfactors in $SU(2)$-invariant Thirring model. Subsequently
we found general formula.
Note that the traces of some concrete operators of such type were
calculated previously in many papers (for example \cite{DFJMN}, \cite{L})
with the help of the Clavelli-Shapiro technique \cite{CS}.

Spaces $\CC [x_1, x_2,...], ~\L [\xi_1,\xi_2,...] $  are called
bosonic and fermionic  Fock spaces, respectively.
They are  the spaces of states in quantum field theory. Observable
values are given
by the  formula  $Tr (\rho A) $, where $\rho$ is the density matrix,
$A$ is the operator of  some observable.
So, we hope that our formulas, will be useful not only
in solving the  integrable models like $SU(2)$-invariant Thirring or
$Sine-Gordon$, but also in the other problems of quantum field theory.

 We shall give several proofs of the formulas \ref{tr_1},\ref{tr_2}.
In the case when  $ A (\partial_{x_1},\partial_{x_2},....), $
$C (x_1, x_2,...) $ are polinomials of the first degree of
 $\partial_1,\partial_2,...$ and $x_1, x_2,...$ respectively,
their commutator equals to the scalar and in this case it  is possible
to prove our formulas in few lines, and this proof is based only
on  commutation relations and the property $Tr AB=Tr BA$.
The idea of another proof applicable in general case is  the
following: if operator $\rho$  is diagonable, then
in the eigenvector basis the sum of diagonal elements
factorizes  to the  product of simple sums.
If the operator is not diagonalizable one can approximate it by
diagonalizable ones, and taking the limit  obtain the formula.
The fact that trace factorizes to the product is based on the following idea:
$\CC [x_1, x_2,...]=\CC [x_1]\otimes\CC [x_2]\otimes... $
and
$\L [\xi_1,\xi_2,...]=\L [\xi_1]\otimes\L [\xi_2]\otimes.... $, and if
the operator $\rho$ is diagonable then
$\rho=\rho_1\otimes\rho_2\otimes... $. And it's well-known that the trace of
the tensor product of operators is product of traces.
The third proof is the longest.  In it we obtain formula for
trace of our operator ristricted on the the space of polinomials of
degree $N$.
Since this space is not tensor product, the sum of diagonal elements
is not factorizable, so the proof is rather intricate.
But lemmas which we proved are rather interesting by themselves.
% and one need them to make rigorous our reasoning about Fredholm's formula.

The alegbras of polinomials $\CC [x_1, x_2,...x_N],\L [\xi_1,\xi_2,...,\xi_N]$
are the synonyms for the symmetrical algebra $ S V$ of linear space $V$
and the exterior algebra $\L V$ of the  space $V$ respectively. The space $V$ is
such that $x_1,x_2,...,x_N$ is a basis of $V$. In this paper it is more
convenient for us to use the  terminology of the symmetrical and exterior algebras
instead of $\CC [x_1, x_2,...x_N],\L [\xi_1,\xi_2,...,\xi_N]$, though, may be it makes
the text not so transparent. For the case of infinite-dimensional $V$
the terminology of the symmetrical and exterior algebras is the only rigorous.

Our formula is valid for the case $dim V=\infty$ as well as $dim V<\infty$.
One can easily reduce infinite-dimensional case to the finite-dimensional one,
because
any operator of the trace class can be approximated by the operators with the
finite-dimensional image.

The paper is organized as follows. In sections 2 and 3  we prove our main
formulas for the cases of fermionic and bosonic Fock spaces, respectively.
Both cases are very similar, so in section 3 we omit all the details.
In the fourth section we obtain Fredholm's formula.
In the fifth section we apply our formulas for calculation of the traces
of vertex operators corresponding to $SU(2)$-invariant Thirring model.
Further on follows several appendices.
In the first we generalize the formulas to the case of the linear space with
a countable basis, without any topology.
In the second we discuss the regularization of the traces.
In the third we mention the case, when there exist a  continuos basis in $V$.
The  fourth appendix contains expressions for the traces of an operators
in Fock spaces in terms of generating functions.
In the fifth appendix we prove that infinite products obtained
in  calculations of the traces of vertex operators in section 5 are
convergent.

{ After the work was completed we received a letter from V.E. Korepin who
pointed out on his works \cite{Kor}, where the Fredholm's determinants
were used to express correlation functions of different integrable
models. This expressions allowed Korepin and his coauthors to obtain
a lot of important information about correlation functions such as
long distance asympotics, equations for correlation functions and so on.}

%%%%%%%%%%%%%%%%%%%%%%%%%%%%%%%%%%%%%%%%%%%%%%%%%%%%%%%%%%%%%%%%%%%%%%%%%%%
%%%%%%%%%%%%%%%%%%%%%%%%%%%%%%%%%%%%%%%%%%%%%%%%%%%%%%%%%%%%%%%%%%%%%%%
%%%%%%%%%%%%%%%%%%%%%%%%%%%%%%%%%%%%%%%%%%%%%%%%%%%%%%%%%%%%%%%%%%%%%%%%%%%%

\setcounter{equation}{0}

\section{Formula for fermionic Fock space.}

In this section we prove the  main trace formula  (theorem 2.1) for
the fermionic Fock space. From the mathematical point of view, fermionic
Fock space is simply the space of polinomials of anticommuting variables or,
more formally, the exterior algebra of some space $V$.
We use the terminology of exterior algebra, because it's more convenient
and it's rigorous in infinite-dimensional case.

 Let $V$ be a linear space. Denote by $\L^n V$
n-th antisymmetrical (exterior) tensor exponent of the space $V$, $\L V=
\oplus_{n=0}^{\infty}\L ^n V$ the exterior algebra of space $V$.
If $a_1, a_2,... $ is a basis of  $V$, then $\L V$ naturally is identified
with  $\L [a_1, a_2,...] $. Any operator $\rho $ on the space $V$
induces the action of operator $\bar\rho $ on the space $\L V$ by the formula:
$\bar\rho (\bigwedge_{i} {\bf v_i})=\bigwedge_i\rho ({\bf v_i}) $, where ${\bf
v_i}$ is  arbitrary elements from $V$. Note that
correspondence $ (V,\rho)\to (\L V,\bar\rho) $ is often called  functor of
secondary quantization \cite{MM}, $\L V$ is called fermionic Fock space.
If $V$ is a Hilbert space, $\rho$ is a unitary operator, then on
completion of $\L V$ one can canonically introduce the structure of Hilbert
space and $\bar\rho$ turns out to be unitary operator on $\L V$. It's natural
to consider algebra $\L V$ graduated: if ${\bf w}\in\L ^n V $ then
$deg {\bf w}=n$.

  Let $V^*$ be a dual linear space. We shall extend the canonical
 pairing between $V^*$
and $V$ to the pairing between $\L V$ and $\L V^*$ according to Vick's rule:
\bea <\bigwedge_{i=1}^m {{\bf f_i}}|\bigwedge_{j=1}^m {{\bf u_j}} >=
\sum_{\sigma\in S_m} (-1) ^{sgn{\sigma}}\prod_{i=1}^{m} < {\bf f_i}| {
\bf u_{\sigma (i) }} > =det|< {\bf f_i}| {\bf u_j} > |_{i, j\leq m}
\label{vik_ext}\label{spar_e_ext}\eea

If $e_i$ and $\tilde e_i$ are dual basis in $V$ and $V^*$, respectively,
then pairing may be defined in the following equivalent way:
$ < (\tilde e_1) ^{i_1}\land...\land (\tilde e_n) ^{i_n}|(x_1) ^{j_1}\land...\land (x_n)
^{j_n} >=\delta_{i_1}^{j_1}...\delta_{i_n}^{j_n} ~~~~~i_k=0, 1~~j_k=0.1 $.

For the arbitrary ${\bf v}\in L V $ one can define
the operator of LEFTWARD multiplication on ${\bf v}$, which we will
call the creating operator and denote it by $C {\bf(v) }$. By definition
operator $C {\bf(v) }$ acts from $\L V$ to $\L V$, as follows:
$C {\bf(v) } [ {\bf \bar v } ]= {\bf v }\land {\bf \bar v}$. Analogically,
for arbitrary ${\bf w}\in L V^* $ one can consider operator
$C({\bf w  })$ : $\L V^* \to \L V^*$. Operator dual to the
$C({\bf w  })$, which acts   $\L V \to \L V$ we will
call the annihilating operator and denote it by $A ({\bf w}) $.
For example, let us describe the action of operator
$A({\tilde e_p})$ on  ${\bf v}=e_{i_1}\land e_{i_2}\land...\land e_{i_l}$, where
$e_i$ is a basis of $V$, $\tilde e_i$ is dual basis of $V^*$ and
$i_k\neq i_l$ at $k\neq l$.
If $p \notin \{i_j\}$ then   operator $A({\tilde e_p})$ kills
${\bf v}$, if $p=i_1$ then  operator erases $e_{i_1}$, if $p \in \{i_j\}$,
but $p\neq i_1$, then one should put $e_p$ on first position, changing
the sign appropriately, then erase it.
If $w\in V^*$, then
$A(w)$ is antiderivating:
$A(w)(v_1\land v_2)=A(w)(v_1)\land v_2+ (-1)^{deg v_1} v_1\land A(w) (v_2)$.
Obviously, the following commutation relations holds:
\bea
&& A ({\bf w}) C ({\bf v}) +C ({\bf v}) A ({\bf w})=< {\bf v}|{
\bf w} > ~~~~~~~~~ {\bf v}\in V, ~~ {\bf w}\in V^*
\label{com_ext}\\
&& A ({\bf w_1}) A ({\bf w_2 })=A ({\bf w_1 w_2}) ~~~~ C ({\bf v_1}) C
({\bf v_2 })=C ({\bf v_1 v_2})
\nn\\&&\bar\rho C ({\bf v_1})=C (\bar\rho ({\bf v_1}))\bar\rho~~~~~ A
({\bf w_1})\bar\rho=\bar\rho A (\bar\rho ({\bf w_1})) ~~~~~~~~~~~~ {
\bf v_1, v_2}\in\L V~~{\bf w_1, w_2}\in\L V^* \nn\\
&&\rho: V\to V,
~~\bar\rho \mbox{ is induced operator} :\L V \to \L V
\nn\eea

Here and further we will allow ourselves
some sloppiness in notations: we will denote
operator dual to $\bar\rho$ by the same symbol $\bar\rho$.
So, in this notations: $ < {\bf w}|\bar\rho{\bf v} >=<\bar\rho{\bf w}|{\bf v} > $.

%We discuss in details infinite-dimensional case
%in appendix annex 1, but for convenience
%we shall formulate here some necessary definitions.

{\bf Definition 2.1:} operator $\rho$ on  $V$
 is called operator with trace, iff
$\exists\l_i\in\CC, ~~y_i\in V, ~~f_i\in V^* $ such that $
\forall x\in V$ holds:

\bea\rho (x)=\sum_{i > 0}\l_i f_i (x) y_i\label{tr_op_def}\eea Where
$\sum_{i > 0}|\l_i| <\infty $ ~~~~ $f_i\in V^*~~~|f_i| < C $ ~~~~ $y_i
\in V~~~ |y_i| < C $

(Small change generalizes this definition to the case of an arbitrary
locally-convex space \cite{SCH}. All results will be true in this case as well.)

{\bf Definition 2.2}: $Tr\rho=\sum_{i > 0}\l_i f_i (y_i)$

Trace is well-defined i.e. trace is independent of the representation \ref{tr_op_def},
if metric Grothendieck's approximation property (AP)
holds  for the space $V$ (\cite{Pietsch}). Let us recall the formulation
of (AP): for any compact subset  $K\subset V$ and $\forall\epsilon\geq 0$
there exists operator $L: |L|\leq 1 $ with finite-dimensional image, such,
that $\forall x\in K ~~~|x-Lx| <\epsilon$.

As it is known, all typical spaces $l^p, L^p[a,b], C^p [a, b] $
are spaces with AP,
the construction  of the space, where AP is not fulfilled,
was long-time standing problem. Known separable examples of such
spaces are rather artificial. Note that if there exists basis (in Schauder
sense) in space $V$, then $V$ is space with AP. We shall also note that
approximation property is closely connected to the following properties:
any compact operator can be approximated in norm topology by operators
with finite dimensional image; identical operator can be approximated
by the compact operators uniformly on any precompact set.

It's clear that if $\rho$ is operator with trace, $A$ is arbitrary operator, then
$\rho A$ and $A\rho$ are operators with trace and $ Tr A\rho= Tr \rho A$.

Here and further on we will always mean that $V$ is a normed space with AP.

{\bf Proposition 2.1} If $\rho$ is a operator
with trace on space $V$, then  $\bar\rho $ is  operator
with trace on $\L V$ and $Tr_{\L V}\bar\rho=\lim{N\to
\infty} Tr_{\L V}\bar\rho_N $ where $\rho_N (x)=\sum_{i > 0}^{N}\l_i
f_i (x) y_i $.

The proof of proposition  is based on straightforward
expression for the $Tr_{\L^k V}\bar\rho$:

\bea
Tr_{\L^k V} \bar \rho = \frac{1}{k!} \sum_{l_1,l_2,...,l_k}
det | f_{l_i} (y_{l_j}) |_{i,j\leq k} \prod_{i=1}^k \l_{l_i}
\label{yvno}
\eea
One can obtain this expression straightforward or using the corollary 2.1 below.
Due to Hadamard's inequality
$det | f_{l_i} (y_{l_j}) |_{i,j\leq k}\leq C^{2k}\sqrt{k}^k$,
hence $ |Tr_{\L^k V} \bar \rho| \leq  \frac{C^{2k} \sqrt{k}^k}{k!} (\sum_i |\l_i|)^k$,
hence due to d'Alembert test the series $\sum_{k} Tr_{\L^k V} \rho $  is
absolutely convergent. Proposition is proved.

{\bf Theorem 2.1}

Let $\rho $ be  operator with trace on $ V$; $ {\bf v}, {\bf w} $
arbitrary elements from $\L V$ and $\L V^*$
respectively. Then the following formulas hold:

\bea
(a) ~~~Tr_{\L V}\left (\bar\rho A ({\bf w}) C ({\bf v})\right)
={Tr_{\L V}\left (\bar\rho\right)}\left < {\bf w}|\frac{1}{1+\bar\rho}
{\bf v}\right >\label{main_f1_ext} \nn \\
(b) ~~~ Tr_{\L V}\left (\bar\rho C ({\bf v}) A ({\bf w})\right)
=Tr_{\L V}\left (\bar\rho\right)\left < {\bf w}|\frac{\bar\rho}{1+\bar
\rho} {\bf v}\right >\label{main_f2_ext}\eea

Remark 1:
In Proposition  3.2 we shall show that
$Tr_{\L V}\bar\rho=det_{V} (1+\rho) $. In accordance with this,
one must understand $Tr_{\L V}\bar\rho\left < {\bf w}|\frac{1}{1+\bar\rho}
{\bf v} \right > $ as matrix element $ < {\bf w}|.... {\bf v} > $
of the augmented  matrix, in the case if $ (1+\rho) $ is not invertable.

Remark 2: Note that the formula \ref{main_f2_ext} (a) is  equivalent
to the formula \ref{tr_2}.

\vskip 2em

{\bf Let us derive \ref{main_f2_ext}$ (b) $ from \ref{main_f2_ext}$ (a) $ }
%{\bf Conclusion\ref{main_f2_ext} (b) from[out of]\ref{main_f1_ext} (a)}

\bea
{Tr_{\L V} \left(
\bar \rho  C({\bf v}) A({\bf w})  \right)}
={Tr_{\L  V} \left(
 A({\bf w}) \bar \rho  C({\bf v})   \right)}
={Tr_{\L  V} \left(
\bar \rho  A(\bar\rho({\bf w})) C({\bf v})   \right)}
=\nn\\
={Tr_{\L  V} \left( \bar \rho \right)} < \bar\rho ({\bf w})| \frac{1}{1+\bar\rho} {\bf v} >
={Tr_{\L  V} \left( \bar \rho \right)} <  {\bf w}| \frac{\bar\rho}{1+\bar\rho} {\bf v} >\nn
\eea

We shall give 3 proofs of the theorem 2.1. The first proof
is very short, %and uses only commutation relations between
%operators $\bar\rho, ~~ A ({\bf w}), ~~ C ({\bf v}) $,
but, infortunately,
it is applicable only for the case $ {\bf w}\in V^*, ~~ {\bf v}\in V $
i.e. $deg {\bf w}=deg{\bf v} =1$. %, in principle one can
%carry out similar proof for arbitrary ${\bf w}, ~~ {
%\bf v}$, but then it will become very bulky.
Another proof is simple too, and applicable in general case.
%  to find the second proof which is simple too.
The third proof is based on three propositions. The proof of the third
is rather complicated, but it is rather interested by itself.
%we need these   propositions not only for the
%proof of the theorem 2.1, but we shall apply them to
%to make our reasoning about Fredholm's formulas rigorous.

Let us formulate the following two simple facts, which are necessary for
the proof of theorem 2.1.

{\bf Lemma 2.1} If theorem 2.1 is true for any operator $\rho$ with trace
such that $|\rho|< 1 $, then theorem 2.1 is true for any $\rho$ with trace.

{\bf Proof:}
Let $\rho$ be an arbitrary operator with trace.
Consider operator $\l\rho$, where $|\l| <\frac{1}{|\rho|}$.
Since $|\l \rho|<1$, we see that theorem 2.1 is true for $\l\rho$, hence:
$ {Tr_{\L V}\left (
\overline{\l\rho} A ({\bf w}) C ({\bf v})\right)}
 = {Tr_{\L V}\left (\overline{\l\rho}\right)} {\left < {\bf w}|
\frac{1}{1+{\l\rho}} {\bf v}\right >} $

Due to proposition 2.1  LHS and RHS of the equality above
exist for all $\l\notin Spec \rho$ and they are analytical functions
of variable $\l$. Hence, due to analicity, equality is true for all $\l\notin Spec \rho$. Taking $\l=1$
one obtains lemma  for arbitrary $\rho$.

\vskip 2em

{\bf Lemma 2.2} If theorem 2.1 is true for arbitrary finite-dimensional
space $V$, then it's true for any space $V$.

{\bf Proof:}
The idea of proof is easy - one can approximate any operator with trace
by operators with finite-dimensional image.

Let $\rho$ be an operator with trace on  $V$.
Let us define operator $\rho_N$ as follows:
 $\rho_N (x)=\sum_{i > 0}^{N}\l_i f_i (x) y_i$,
$\l_i,~y_i,~,f_i$ are defined according to definition 2.1.
Let space  $ V_N$ be
linear span of vectors $v, y_1, y_2,... y_N$. Then $Im\rho_N\subset
V_N$ that is why $Tr_{\L V}\rho_N=Tr_{\L V_N}\rho_N$. Since formula
is true for any finite-dimensional space, hence
$Tr_{\L V_N}\left (\bar\rho_N A
({\bf w}) C ({\bf v})\right)= {Tr_{\L V_N}\left (\bar\rho_N\right)}
\left < {\bf w}|\frac{1}{1 -\bar\rho_N} {\bf v}\right >$

Hence:
\bea
\frac{Tr_{\L V}\left (\bar\rho_N A ({\bf w}) C ({\bf v})\right)
}{Tr_{\L V}\left (\bar\rho_N\right) }=\frac{Tr_{\L V_N}\left (\bar\rho_N A
({\bf w}) C ({\bf v})\right) }{Tr_{\L V_N}\left (\bar\rho_N\right)}
=\left < {\bf w}|\frac{1}{1 + \bar\rho_N} {\bf v}\right >\nn
\eea
Taking the limit $N\to\infty$ one obtains the necessary lemma.
Passage to the limit is valid, because of proposition 2.1.

{\bf Lemma 2.2 is proved.}

{\bf The first Proof of theorem 2.1.}  Only for case
$deg {\bf w}=deg {\bf v}=1 $. ( It is sufficient for Fredholm's formulas.)

Recall formula \ref{com_ext}: $A ({\bf w}) C ({\bf v}) +C ({
\bf v}) A ({\bf w})=< {\bf v}|{\bf w} > $

Due to lemma 2.1 it's sufficient to prove
only for such $\rho$ that $\rho^n\to 0,~~ n\to\infty $.

\bea &&Tr_{\L V}\left (\bar\rho A ({\bf w}) C ({\bf v})\right)
=\left < {\bf w} |{\bf v}\right > Tr_{\L V} (\bar\rho)-
Tr_{\L V}\left (\bar\rho C ({\bf v}) A ({\bf w})\right)\nn\\&&
={\left < {\bf w} |{\bf v}\right >} Tr_{\L V} (\bar\rho)-
Tr_{\L V}\left (C (\rho ({\bf v}))\bar\rho A ({\bf w})\right)
={\left < {\bf w} |{\bf v}\right >} Tr_{\L V} (\bar\rho)-
Tr_{\L V}\left (\bar\rho A ({\bf w}) C ({\rho (\bf v) })\right)\nn\\
&&\mbox{ proceeding this transformation n times one gets:}\nn\\
&&Tr_{
\L V} (\bar\rho)\sum_{i=0}^{n-1} {\left < {\bf w}|(-\rho) ^i {\bf v}
\right > }+ (-1) ^nTr_{\L V}\left (\bar\rho A ({\bf w}) C ({\rho^n (
\bf v) })\right)
=Tr_{\L V} (\bar\rho) {\left < {\bf w}|\frac{1}{1+\rho} {\bf v}\right >}
\nn\eea
The latest equality is obtained
by sending $n\to\infty$, which is valid since $\rho^n\to 0 $.

{\bf\mbox{Proof 1 is finished.}}

Remark 1:
The difficulty of applying this proof to the case, when $deg {\bf w} >1 , deg {\bf v}>1 $,
is the following.
In this case commutation relations between $A({\bf w}),C({\bf v})$ are complicated.
So the way to prove  described  above is possible, but very bulky,
and we prefer to give another proof.

Remark 2: The proof that traces of vertex operators satisfy
Kniznik-Zamolodchikov equations is similar to proof 1 described above.

Remark 3: Let us note that we have not used that $dim V < \infty$.
We have only used commutation relations between
$\bar\rho, ~A ({\bf w}), ~ C ({\bf v}) $ and  the property that $TrAB=TrBA $.
But this is not surprising, because  trace maybe
determined by the property that $TrAB=TrBA $.

{\bf The second proof.}

Due to  lemma 2.2 it's sufficient to prove only for $dim V=N <\infty$.

Proof is based on two simple ideas. The first idea is the following:
it's sufficient to consider the case when $\rho$ is diagonable operator,
because any not diagonable operator can be approximated by diagonable ones.
And easy to see that if $\rho_i\to \rho$ and theorem holds for $\rho_i$,
then it holds for $\rho$. (Due to $dimV<\infty$, hence
$dim \L V< \infty$, we have no problems with approximation).
The second idea is the following: if one writes down the sum of diagonal
elements in the basis of eigenvectors for $\rho$, then this sum factorizes. And
each multiplier is the series similar to the geometric progression.

Let us consider the case $dimV=1$. Let $e\neq 0\in V$ ; ~~ $\rho e=\l e$,
let $\tilde e$ be the dual to $e$ element of $V^*$

\bea &&Tr_{\L V} (\bar\rho)=< 1|\rho^0 1 >+<\tilde e|\rho e > =1+\l
=det_{V} (1+\rho)\nn\\&& Tr_{\L V} (\bar\rho A (\tilde e) C (e))=< 1|
\rho^0\partial_e (e\land 1) >+<\tilde e|\rho\partial_e (e\land e) >=
\nn\\&&
=1+0 =1= (1+\l) <\tilde e|\frac{1}{1+\l}e >=Tr_{\L V}\bar\rho <\tilde
e|\frac{1}{1+\rho}e >\label{lamb}\\&& Tr_{\L V} (\bar\rho A (\tilde e)
C (1)) =Tr_{\L V} (\bar\rho A (1) C (e))=<\tilde e|\frac{1}{1+\rho} 1
>=< 1|\frac{1}{1+\rho} e > =0\nn\\
&&\mbox {Since}<\L^k V|\L^i V > =0 \mbox { for } i\neq k\nn
\eea

In proving : we have obtained that
$Tr_{\L V}\bar\rho={1 +\l}$. Note that one can
rewrite it  as follows: $Tr_{\L V}\bar\rho =det_{V}(1 + \rho)$.
Further  we shall prove that it's true not only for $dimV=1$, but also
for any $V$ (even infinite-dimensional).

We see that transformations in formula \ref{lamb} are possible
only in the case $\l\neq-1$, i.e. $1+\rho$ is invertable.
Otherwise one must stop transformations on the step:
$ Tr_{\L V} (\bar\rho A (\tilde e) C (e)) =1$, and
understand it as  matrix element $ <\tilde e|...e > $ of
matrix augmented to $1+\rho$, easily to see that this
will be true in case $dim V > 1$.

Thus theorem 2.1 is proved for the case $dim V=1$.

\vskip 1em

Let us  come back to the case   $dimV=N $ where N is arbitrary.

Recall that we want to prove:
\bea
Tr_{\L V} \left( \bar \rho {A({\bf w})} {C({\bf v})} \right)
= \left(Tr_{\L V}  \bar \rho \right)
\left< {\bf w}~| \frac{1}{1+ \bar \rho} {\bf v} \right>\nn
%\label{fml2}
\eea

Note that LHS and RHS of the formula above are bilinear
of  ${\bf w}, {\bf v}$. Therefore it's sufficient
to prove only for basic vectors.

Let $c_1, c_2,..., c_N$ be the basis of eigenvectors for operator $\rho$,
(as we have already  said it's sufficient
to prove for diagonable opeartors only).
Let $\l_i$ be the eigenvalues of operator $\rho$ corresponding to
vectors $c_i$. Let $\tilde c_1,\tilde c_2,...,\tilde c_N$ be the dual
basis in $V^*$.

Exterior algebra of N variables is the tensor product
of the exterior algebras of each variable,
therefore
$\L V=\L [c_1]\ot\L[c_2]\ot....\ot\L [c_N]$. Easy to see that
$\bar \rho
\prod_{i=1}^N {A(\tilde c_i)}^{r_i} \prod_{j=1}^N {C(c_j)}^{t_j}
=\bigotimes_{i=1}^N \bar\rho_i
{A(\tilde c_i)}^{r_i} {C(c_i)}^{t_i}
$
where
$\bar \rho_i$ acts in  $\L[e_i]$  as follows: $\bar\rho_i(e_i^k)=\bar\rho(e_i^k)=\l_i^k (e_i^k)$
~~ and ~~~($r_i,t_j, k = 0,1 $)

That is why:
\bea
&&Tr_{\L V}\left (\bar\rho\prod_{i=1}^N {A (\tilde c_i) }^{r_i}
\prod_{j=1}^N {C (c_j) }^{t_j}\right)=\prod_{i=1}^{N} Tr_{\L [c_i]}
\left (\bar\rho_i {A (\tilde c_i) }^{r_i} {C (c_i) }^{t_i}\right)=\nn\\
&&=\prod_{i=1}^{N} (Tr_{\L [c_i]}\bar\rho_i) < {(\tilde c_i) }^{r_i}|
\frac{1}{1+\bar\rho_i}
={(c_i) }^{t_i} >  (Tr_{\L V}\bar\rho) <\prod_{i=1}^{N} {(\tilde c_i)
}^{r_i}|\frac{1}{1+\bar\rho}\prod_{j=1}^{N} {({c_j}) }^{t_j} >\nn\\
&&\mbox{The last equality is true due to: }
\rho\mbox{ diagonable in basis } c_i\mbox{ and} < c_i|c_j> =0 ~i\neq j\nn
\eea

{\bf Proof 2 is finished.}

Note that
 $Tr_{\L V}\bar\rho=\prod_{i=1}^N Tr_{\L [c_i]}\bar\rho_i=\prod_{i=1}^N (1+\l_i)=det_{V} (1+\rho) $

Thus we have proved the following, well known proposition:

{\bf Proposition 2.2} Let $V$ be a finite-dimensional space, let
$\rho$ be an operator on $V$, then
\bea det_{V} (1+\rho) =Tr_{\L V}\bar\rho\label{det}\eea

There is another more direct proof of this proposition:

Let  dim V=N, then by definition of determinant
$det_V O=Tr_{\L^N V}\bar O$. Hence:
$det_V (1+\rho) =Tr_{\L^N} (1+\rho) ^{\land n}$
$=Tr_{\L^0}\bar\rho
+Tr_{\L^1}\bar\rho +Tr_{\L^2}\bar\rho+..... +Tr_{\L^N}\bar\rho$

%\bea &&\mbox{Let  dim V=N, then by definition of determinant }
%det_V O=Tr_{\L^N V}\bar O\nn\\&&\mbox{Hence:}
%det_V (1+\rho) =Tr_{\L^N} (1+\rho) ^{\land n}=\nn\\
%%&&=Tr_{\L^n} (1\land 1\land...\land 1)+Tr_{\L^n}\left (\rho\land 1
%%\land 1\land... 1+1\land\rho\land 1\land 1...\land 1+.... +1\land 1
%%\land...\land\rho\right)+\nn\\&&+Tr_{\L^n}\left (\rho\land\rho\land 1
%%\land 1...\land 1+\rho\land 1\land\rho\land 1.. 1+.... +1\land 1\land 1...
%%\land\rho\land\rho\right)+\nn\\&&+.............. +Tr_{\L^n} (\rho
%%\land\rho\land\rho\land\rho\land\rho...\land\rho)
%&&=Tr_{\L^0}\bar\rho
%+Tr_{\L^1}\bar\rho +Tr_{\L^2}\bar\rho+..... +Tr_{\L^N}\bar\rho\nn\eea

%So  formula $det_V (1+\rho)=\sum_{n=0}^{\infty}Tr_{\Lambda^n V}\rho$
%is the direct consequence of the definition.
%The author doesn't know if there is any analogous direct proof
%of proposition 2.1: $det_V\frac{1}{1 -\rho}=\sum_{n=0}^{\infty}Tr_{S^n V}\rho$.

Motivated by  \ref{det}  one can introduce the following
definition of determinant (applicable in infinite-dimensional case):

{\bf Definition 2.3} Operator $O$ on space $V$ is called by
operator with determinant, if operator $\rho=O-1$ is  operator with
trace. Determinant of  $O$ is defined by the formula: $ det_{V} O=Tr_{\L V}\bar\rho$.

Further we shall show
%that this definition conincides with definition of Fredholm's determinant.
that this definition is  Fredholm's definition, but
written in more invariant way.
Note that this is definition of usual, not somehow regularized, determinant.
This determinant is equal to the product of the eigenvalues of operator $O$.

\vskip 1em

{\bf The third proof.}

{\bf Proposition 2.3}: Let $e_i$ be  an arbitrary basis of
the finite-dimensional space $V$, $\tilde e_i$ be the dual basis,
then for any operator $O$ on space $\L^n V$,
the following formula takes place:

\bea Tr_{\Lambda^n V} O=\frac{1}{n!}\sum_{l_1... l_n\geq 0} {<\tilde
e_{l_1}\land\tilde e_{l_2}\land...\land\tilde e_{l_n}|O e_{l_1}\land
e_{l_2}\land...\land e_{l_n} > }&&\label{sim_sim_ext}\eea

The proof is obvious, due to each basic
vector is repeated in summation $n!$ times and  it cancels with
$n!$ in denominator.

Remark: If $V$ is Hilbert space, $e_i$ is it's basis, $O$ is operator
with trace on $\L^n V$, then
lemma is true also. The proof is the same, but instead of finite
sums arises series, which are absolutely convergent, since $O$ is
operator with trace.

{\bf Corollary 2.1:}  Let $e_i$ be  an arbitrary basis of
the finite-dimensional space $V$, $\tilde e_i$ be the dual basis.
 Let $\rho$ be an operator
on space $V$, $\bar\rho$ - induced operator on $\L^n V$. Then:

\bea Tr_{\Lambda^n V}\bar\rho=\frac{1}{n!}\sum_{l_1... l_n\geq 0} det
\left\vert\matrix{ <\tilde e_{l_1}|\rho e_{l_1} > & <\tilde e_{l_2}|
\rho e_{l_1} > &... & <\tilde e_{l_n}|\rho e_{l_1} >\cr <\tilde e_{l_1}|
\rho e_{l_2} > & <\tilde e_{l_2}|\rho e_{l_2} > &... & <\tilde e_{l_n}|
\rho e_{l_2} >\cr
\vdots &\vdots &\ddots &\vdots\cr <\tilde e_{l_1}|\rho e_{l_n} > & <
\tilde e_{l_2}|\rho e_{l_n} > &... & <\tilde e_{l_n}|\rho e_{l_n} >}
\right\vert &&\label{sim_sim_det}\eea
To prove corollary one should apply Vick formula \ref{vik_ext} to \ref{sim_sim_ext}.

{\bf Proposition 2.4} Let $\rho$ be an operator with trace on space
 $V$, $\bar\rho$ induced operator on $\L^n V$,
then $Tr_{\L^n V}\bar\rho $ can be expressed
through  $Tr_{V}\rho^l $ as follows:

\bea Tr_{\L^n V}\bar\rho=\sum_{ n_j\geq0:\sum_j jn_j=n}
\frac{  (-1)^{\sum j(n_j-1)}}{n_1! n_2!... 1^{n_1}2^{n_2}...}\prod_{l} (Tr_V (\rho^l)) ^{n_l}\eea

{\bf  Proposition 2.5} Let $\rho $ be an  operator with trace on $ V$;
$ {\bf v_i}, {\bf w_i} $ arbitrary elements from  $V$ and $V^*$,
respectively. Then the   trace
of $\bar\rho\prod_{i=1}^k {A ({\bf w_i})}\prod_{j=1}^k {C ({\bf v_j}) }$
over $\L^nV$ can be expressed through $Tr_{\L^r V}\bar\rho$ and $ < {\bf w_i}|~\bar\rho ^{l}~{\bf v_j} > $
as follows:

\bea Tr_{\L^n V}\left (\bar\rho\prod_{i=1}^k {A ({\bf w_i})}
\prod_{j=1}^k {C ({\bf v_j})}\right)=\sum_{\begin{array}{c} r, q_j
\geq0:\\\sum_j q_j+r=n\end{array}} Tr_{\L^r V} (\bar\rho)\left <
\prod_{j=1}^k {\bf w_j}|\prod_{j=1}^k (-\rho) ^{q_j} ({\bf v_j})\right
> &&\label{lem3_ext}\eea

{\bf Proof of theorem 2.1 easily follows from Proposition 2.5:}

LHS and RHS of equality in theorem 2.1 are bilinear of {\bf w}, {\bf v},
therefore it's sufficient to prove only  for
${\bf w}=\bigwedge_{i=1}^k {{\bf w_i}}, ~~~{\bf v}=\bigwedge_{j=1}^n {{\bf v_j}}$,
where $deg {\bf w_i}=deg {\bf v_i}=1$. Easily to see that
both LHS and RHS are equal to zero, if $n\neq k$.
So it's sufficient to consider the case $n=k$ only.

Due to lemma 2.1 it's sufficient to prove for such $\rho$ that $|\rho|<1$
\bea && Tr_{\L  V}\left (\bar\rho\prod_{i=1}^k {A ({\bf w_i})}
\prod_{j=1}^k {C ({\bf v_j})}\right)
=\sum_{n\geq 0} Tr_{\L ^n V}\left (\bar\rho\prod_{i=1}^k {A ({\bf w_i})}
\prod_{j=1}^k {C ({\bf v_j})}\right)\nn\\
&&=\sum_{n\geq 0}\sum_{ r,n_1,... n_k\geq0: ~\sum_j n_j+r=n}
\left ((Tr_{\L ^rV}\bar\rho) <\prod_{p=1}^k{\bf w_p }~|~\prod_{p=1}^k
(-\rho)^{n_p} ({\bf v_p}) >\right)
\nn\\&&
=\sum_{n_1,..., n_k, r\geq0}\left ((Tr_{\L ^rV}\bar\rho) <\prod_{p=1}^k {
\bf w_p }|\prod_{p=1}^k(-\rho)^{n_p} ({\bf v_p}) >
\right)\nn\\
&&=\left (\sum_{r\geq0} Tr_{\L ^r V}\bar\rho\right)\left <
\prod_{p=1}^k {\bf w_p }|\prod_{p=1}^k\sum_{n\geq 0}(-\rho)^{n} ({\bf v_p})
\right >
=\left (Tr_{\L  V}\bar\rho\right)\left <\prod_{p=1}^k {\bf w_p }|
\frac{1}{1+\bar\rho}\prod_{p=1}^k {\bf v_p}\right >\nn
\eea
{\bf The third proof of the theorem 2.1 is finished.}

{\bf The proof of Proposition 2.4}
Similar to the proof of lemma 2.2, easy to show that it's
sufficient to consider finite-dimensional spaces $V$ only, because
any operator with trace can be approximated by the operators with
finite-dimensional image.
\bea
 &&Tr_{\L ^n V}\bar\rho=
\frac{1}{n!}\sum_{l_1... l_n\geq 0} {<\tilde e_{l_1}\land \tilde e_{l_2}\land...
\tilde e_{l_n}|\bar\rho e_{l_1}\land e_{l_2}\land... e_{l_n} > }=
\nn\\&&=
\frac{1}{n!}\sum_{l_1... l_n\geq 0}\sum_{\sigma\in \L _n}
(-1)^{sgn(\s)}\prod_{i=1}^{n} < e_{l_i}|\rho e_{\sigma (l_i)} >\nn\\&&
=\frac{1}{n!}\sum_{\sigma\in \L _n} (-1)^{\sum i(k_i-1)}
\prod_{i} (Tr_{V}\rho^i) ^{k_i (\sigma)}
\nn\\&&
\mbox{ where } k_i (\sigma)\mbox{ is the number of the cycles of length i for }\sigma\nn\\&&
\mbox{ The amount of permutations that are product of } k_i
\mbox{ cycles of length i, equals to }\frac{n! }{k_1! k_2!... 1^{k_1}2^{k_2}...}\nn\\
&&\mbox{ hence, equality may be continued:}
\nn\\
&&=\sum_{k_j\geq0:\sum_j jk_j=n}
\frac{ (-1)^{\sum i(k_i-1)}}{k_1! k_2!... 1^{k_1}2^{k_2}...}\prod_{i} (Tr_V\rho^i) ^{k_i}
\nn\eea
{\bf Proposition 2.4 is proved}

{\bf The proof of Proposition 2.5}
Similar to the proof of lemma 2.2, easy to show that it's
sufficient to consider only finite-dimensional spaces $V$ only, because
any operator with trace can be approximated by the operators with
finite-dimensional image.

The main instrument of proof is  - is Vick's formula (\ref{vik_ext}).
Note that Vick's formula is conveniently  rewritten rule
of pairing between antisymmetrical tensors.

\bea && Tr_{\L^n V}\left (\bar\rho\prod_{i=1}^k {A ({\bf w_i})}
\prod_{j=1}^k {C ({\bf v_j})}\right)
=\nn\\
&&=\frac{1}{n!}\sum_{l_1... l_n\geq 0} {<\tilde e_{l_1}\land \tilde e_{l_2}\land...
\land\tilde e_{l_n}|
\bar\rho\prod_{i=1}^k {A ({\bf w_i})}\prod_{j=1}^k {C ({\bf v_j})}
e_{l_1} \land e_{l_2}\land...\land e_{l_n} > }=\nn\\&&
=\frac{1}{n!}\sum_{l_1... l_n\geq 0} {< \bigwedge_{i=1}^k {\bf w_i}\land
\rho (\tilde e_{l_1})\land\rho (
\tilde e_{l_2})\land...\land\rho (\tilde e_{l_n})|
\bigwedge_{j=1}^k {\bf v_j}\land  e_{l_1}\land  e_{l_2}\land ... \land e_{l_n} > }=\label{n_1}
\eea

Let us introduce the following notation:

\bea
&&
{\xi_{l_i}}={\bf v_{i}} ~~, ~~
{ \tilde  \xi_{l_i}}= {\bf w_{i}} ~~, ~~1 \leq i \leq k ;
~~~~~~~~~~
{ \xi_{l_{i}}}=e_{l_{i+k}} ~~ , ~~{ \tilde  \xi_{l_{i}}}
=\rho (\tilde e_{l_{i+k}})~~ , ~~k < i \leq n+k ;\nn\\
&&\mbox{Let us apply Vick's formula to the \ref{n_1}:}\nn\\
&&=\frac{1}{n!}\sum_{l_1... l_n\geq 0}\sum_{\sigma\in S_{n+k}} (-1)^{sgn \s}
\prod_{i=1}^{n+k} < {\bf\tilde\xi_{l_i}}|{\bf\xi_{l_{\sigma (i) }}} >
\nn\\
&&\mbox{Rewrite multipliers in following order:}\nn\\
&&=\frac{1}{n!}
\sum_{l_1...l_n \geq 0} \sum_{\sigma \in S_{n+k}} (-1)^{sgn \s} \nn\\
&&
<{\bf w_1 }|  e_{\bf{l_{\sigma(1)}}} >
< \rho(\tilde e_{\bf{l_{\sigma(1) }}}  ) | e_{\bf{l_{\sigma^2(1)}}} > ....
< \rho( \tilde e_{\bf {l_{\sigma^{n_1}(1)  }}} ) | {\bf v_{\sigma^{n_1+1}(1)}}  > \nn\\
&&<{\bf w_2 }|e_{\bf {l_{\sigma(2)}}} >
< \rho(\tilde e_{\bf{l_{\sigma(2)}}})| e_{\bf{l_{\sigma^2(2)}}} > ....
< \rho(\tilde e_{\bf {l_{\sigma^{n_2}(2)}}})| {\bf v_{\sigma^{n_2+1}(2)}}  > \nn\\
&&.............................. \nn\\
&&<{\bf w_k }|  e_{\bf{l_{\sigma(k)}}} >
< \rho(\tilde e_{\bf{l_{\sigma(k)}}} ) | e_{\bf{l_{\sigma^2(k)}}} > ....
< \rho(\tilde e_{\bf {l_{\sigma^{n_k}(k)}}})| {\bf v_{\sigma^{n_k+1}(k)}}  > \nn\\
%{\bf \xi_{l_{\sigma(i)}}}>\nn\\
&&\prod_{j\in {\bf J(\sigma)} } <\rho(\tilde e_{\bf l_j} ) |e_{\bf l_{\sigma(j)}}>
\eea

Where ${\bf J (\s) }$ is subset of \{1... n\} that includes such elements
that cannot be obtained
from \{1, 2,..., k\} by applying  $\sigma,\s^2,\s^3... $.

$r=Card (\J) $.~~~~~~~~( As usual, Card(M) is amount of elements in the set M ).

$n_i\geq 0$ - the least  exponent such that  $\s^{n_i+1} (i)
$ belongs to the set $\{1, 2,3,... k\}$.

Denote ${\bf I (\s) }$=\{1, 2,... n\} $\backslash$ ${\bf J (\s) }$.

$Card (\I)=\sum_{i} n_i$.

Then $r+\sum_i n_i =n$.

Summating over $l_i$ such that $i\in {\bf I (\s) }$ one obtains:
\bea
=\frac{1}{n!}\sum_{\sigma\in S_{n+k}} (-1)^{sgn \s}
\prod_{p=1}^k < {\bf w_p }|\bar\rho^{n_p} {\bf v_{
\sigma^{n_p+1} (p)}} >
\sum_{~~l_j, j\in
\J}\prod_{j\in\J} <\rho (\tilde e_{\bf l_j}) |e_{\bf l_{\sigma (j) }} >
\nn\eea

Let $\s=\s_1\s_2=\s_2\s_1$, where $\s_1$ is the product of
independent cycles which includes \{1, 2,..., k\}, $\s_2$ -
the one which doesn't includes \{1, 2,..., k\}. Obviously
$\s_1 \s_2=\s_2 \s_1$ and $\s(i)=\s_1(i)~~~i \in \I~~~~\s(i)=\s_2(i)~~~i \in \J$.

Then
\bea &&=\frac{1}{n!}\sum_{\s_1} (-1)^{sgn \s_1}
\left(\prod_{p=1}^k < {\bf w_p }|
\rho^{n_p} {\bf v_{\sigma^{n_p+1} (p)}} >\sum_{~~l_j, j\in\J}
\sum_{\s_2} (-1)^{sgn \s_2}
\prod_{j\in\J} <\tilde e_{\bf l_j}|\rho e_{\bf l_{\sigma (
j) }} >\right)\nn\\&&=\frac{1}{n!}\sum_{\s_1} (-1)^{sgn \s_1}
\left (\prod_{p=1}^k < {
\bf w_p }|\rho^{n_p} {\bf v_{\sigma^{n_p+1} (p)}} > r! Tr_{\L^rV}
\bar\rho\right)\nn\\&&=\frac{r! }{n!}\sum_{\tilde\s\in S_k}
(-1)^{sgn \tilde \s}
\sum_{r,n_1... n_k\geq0:\sum_j n_j+r=n}
\sum_{\s_1\in\Omega (r, n_i,\tilde\s)}\left (\prod_{p=1}^k < {\bf w_p
}|(-\rho)^{n_p} {\bf v_{\sigma^{n_p+1} (p)}} > Tr_{\L^rV}\bar\rho
\right)\nn\\&&=\frac{r! }{n!}\sum_{\tilde\s\in S_k} (-1)^{sgn \tilde \s}
\sum_{r, n_1... n_k\geq0:\sum_j n_j+r=n}
\frac{n! }{r!}\left (\prod_{p=1}^k < {\bf w_p }|(-\rho)^{n_p} {\bf v_{
\sigma^{n_p+1} (p)}} > Tr_{\L^rV}\bar\rho\right)\nn\\&&=\sum_{
\tilde\s\in S_k} (-1)^{sgn \tilde \s}
\sum_{r, n_1... n_k\geq0:\sum_j n_j+r=n}
\left (\prod_{p=1}^k < {\bf w_p }|(-\rho)^{n_p} {\bf v_{\tilde\s (p) }} >
Tr_{\L^rV}\bar\rho\right)\nn\\&&
=\sum_{r, n_1... n_k\geq0:\sum_j n_j+r=n}
\left (Tr_{\L^rV}\bar\rho <\bigwedge_{p=1}^k {\bf w_p }|\bigwedge_{p=1}^k
(-\rho)^{n_p} ({\bf v_p}) >\right)\nn\eea

Where $\Omega (r, n_i,\tilde\s) $ is set of  that and only that
permutations $\s_1$ from  $S_{k}$, such that: $\s_1^{n_i+1} (i) =\tilde\s (
i) ~~; ~~\s_1 (j) =j ~~j\in {\bf J (\s_1) }$, and $n_i $ - are the least
exponents such that $\s_1^{n_i+1}(i)=\tilde\s(i) $
belongs to the set \{1, 2,..., k\}.

Easy to see that $Card (\Omega (r, n_i,\tilde\s))=\frac{n!}{r! }$

{\bf Proposition 3 is proved.}

Due to Corollary 2.1 one can reformulate definition 2.3 as follows:

{\bf Proposition 2.6:} Let $\rho$ be an operator with trace on
finite-dimensional space $V$, $e_i$ a basis of $V$,$\tilde e_i$ be the
dual basis.

\bea
det_{V} (1+\rho)=\sum_{n=0}^{\infty}\frac{1}{n!}\sum_{l_1, l_2,...,
l_n} det|<\tilde e_{l_j}|\rho e_{l_i} > |_{i, j\leq n}
\label{det_2}
\eea

%\section{Formula for Boson Fock.}%\section{Formula for Boson Fock.}
%\section{Formula for Boson Fock.}%\section{Formula for Boson Fock.}
%\section{Formula for Boson Fock.}%\section{Formula for Boson Fock.}

\setcounter{equation}{0}
\section{Formula for bosonic Fock space.}

In this section we prove the  main trace formula  (theorem 3.1) for
the bosonic Fock space. From the mathematical point of view, bosonic
Fock space is simply the space of polinomials or,
more formally, the symmetrical algebra of some space $V$.
We use the terminology of the symmetrical algebra, because it's more convenient
and it's rigorous in infinite-dimensional case.

 Let $V$ be  a linear space. % Denote by $S^n V$ - n-th
$V^{\otimes i}$ - i-th tensor power of the space V, $S^iV$ - subspace
of symmetric tensors,
$SV=\oplus_{i=0}^{\infty} S^i V$ - symmetrical algebra of space V,
If space $V$ is finitedimensional and $e_1, e_2,..., e_N$ - its basis,
then $SV=\CC [e_1, e_2,..., e_N] $.
If $V$ is infinitedimensional then one may naively treat $SV$ as
$\CC[[e_1,e_2,...,]$.
Any operator $\rho $ on the space $V$ induce the action of operator
$\bar\rho $ on the space $S V$ by the  formula:
$\bar\rho (\prod_{i} {\bf v_i})=\prod_i\rho ({\bf v_i}) $ where $v_i$ -
arbitrary elements from $V$. Note that correspondence
$ (V,\rho)\to (SV,\bar\rho) $  is often called  secondary quantization
functor \cite{MM}, $SV$ - bosonic Fock space and if $V$ -
Hilbert space and $\rho$ unitary operator, then on
completion of $SV$ one can canonically introduce the structure
of Hilbert space and $\bar\rho$ turns out to be
unitary operator on it. Naturally to consider
algebra $S V$ graduated: if ${\bf w}\in S^n V $ that $deg {\bf w}=n$.

%%%%Ot nachala seccii do suda peredelano za 50 minut.

  Let $V^*$ be  the dual space to $V$. Extend the  pairing between $V$ and $V^*$
to the pairing between
$S V$ and $S V^*$ according to Vick's rule:

\bea <\prod_{i=1}^m {{\bf f_i}}|\prod_{j=1}^m {{\bf u_j}} >=
\sum_{\sigma\in S_m}\prod_{i=1}^{m} < {\bf f_i}| {\bf u_{\sigma (i) }}
>=per |< {\bf f_i}| {\bf v_j} > |_{i,j\leq m}
~~~~~~~~~~~~ {\bf f_i}\in V^*, ~~{\bf u_j}\in V\label{vik}
\label{spar_e}\eea
where $per A$ is permanent of matrix A. Recall that permanent of square
matrix A is sum of products of its elements, standing in different rows and
columns.

Let us point out that $<{\bf w}^k|{\bf v}^k>=k!<{\bf w}|{\bf v}>^k$,
at  ${\bf w} \in V^*$, ~~~ ${\bf v} \in V$.  It's  convenient, because,
if $x_i$ is basis in $V$, $\partial_i$ is dual basis in $V^*$ then
differential operator $\frac{\partial}{\partial {x_j}}$ in space $SV$
is dual to operator of multiplication on $\partial_{j}$ in space $S V^*$, i.e.:

$< (\partial_1)^{i_1}...(\partial_k)^{i_k+1}...(\partial_n)^{i_n}~~ | (x_1)^{j_1}...(x_n)^{j_n} >=
< (\partial_1)^{i_1}...(\partial_k)^{i_k}...(\partial_n)^{i_n}~|~ \frac{\partial}{\partial {x_k}} ~[ (x_1)^{j_1}...(x_n)^{j_n}] >$

Analogically operator of multiplication on $x_j$ in space $SV$ is dual to differential
operator with respect to $\partial_j$ in space $SV^*$. Also, due to such choice of pairing
proposition 3.3 holds.

Similar as it was done in previous section one can introduce creating-annihilating
operators. For any ${\bf v}\in SV$ creating operator $C({\bf v})$ : $SV\to SV$,
for any ${\bf w}\in SV^*$ annihilating operator $A({\bf w})$ : $SV\to SV$.
In terms of space of polinomials creating operator is operator of multiplication
on some polinomial, annihilating operator is polinomial of differential operators
$\frac{\partial}{\partial x_j}$. Obviously, the following commutation relations holds:
\bea
&&{ A({\bf w})C({\bf v})-C({\bf v}) A({\bf w}) =<v|w> } ~~~~~~~~~ {\bf v}\in V , ~~ {\bf w}\in V^*
\label{com}\\
&& A({\bf  w_1})A({\bf w_2 })=A({\bf w_2}) A({\bf w_1})= A({\bf w_1 w_2})
~~~~ C({\bf  v_1})C({\bf v_2 })=C({\bf v_2}) C({\bf v_1})= C({\bf v_1 v_2})
\nn\\
&&
\bar\rho C({\bf v_1})= C(\bar\rho({\bf v_1})) \bar\rho~~~~~
A({\bf w_1}) \bar\rho= \bar\rho A(\bar\rho({\bf w_1}))
~~~~~~~ {\bf v_1,v_2} \in SV~~{\bf w_1,w_2} \in SV^* \nn\\
&&\rho: V\to V,
~~\bar\rho \mbox{ is induced operator} :SV \to SV
\nn
\eea

{\bf Proposition 3.1} If $\rho$ is a operator
with trace on space $V$ and $|\rho|<1 $, then  $\bar\rho $ is  operator
with trace on $S V$. For $ {\bf v}\in S V , {\bf w}\in S V^* $
operator $\bar\rho A({\bf w}) C({\bf v})$ is operator with trace.

Let us point out that analogous proposition 2.1 in previous section
is true for all operators $\rho$ with trace, in contrast to the proposition 3.1,
which is true only for $|\rho|<1 $. If $|\rho|>1 $, then $\bar\rho$ is
unbounded operator on $SV$. Also note that norm of operator
$A({\bf w})$ on space $\oplus_{i<N} S^i V$ equals to $N |{\bf w}|$.
so it is unbounded on the hole space $S V$, but product
$\bar\rho A({\bf w}) $ is operator with trace.

One can easy derive the proof of the proposition 3.1 from proposition 3.4.

{\bf Theorem 3.1}

Let $\rho $ be an  operator with
trace on $ V$ and $|\rho|<1$ ; $ {\bf v}, {\bf w} $ arbitrary elements from
$SV$ and $SV^*$ respectively.
Then the following formulas take place:
\bea (a) ~~~~~~\frac{Tr_{SV}\left (\bar\rho A ({\bf w}) C ({\bf v})
\right) }{Tr_{SV}\left (\bar\rho\right)}
=\sum_{n=0}^{\infty}\left < {\bf w}|\bar\rho^n {\bf v}\right >
~~~~~~~~~~ (b) ~~~~~~\frac{Tr_{SV}\left (\bar\rho C ({\bf v}) A ({\bf
w})\right) }{Tr_{SV}\left (\bar\rho\right)}
=\sum_{n=1}^{\infty}\left < {\bf w}|\bar\rho^n {\bf v}\right >
\label{main_f2}\label{main_f1}\eea
Remark 1: Obviously, formula \ref{main_f2} (a) is equivalent to the formula
 \ref{tr_1}, mentioned in introduction.

Remark 2: Actually if exist two of the expressions $ Tr_{SV} (\bar\rho), $
 $Tr_{SV}\left (\bar\rho A ({\bf w}) C ({\bf v})\right), $ $
\sum_{n=0}^{\infty}\left < {\bf w}|\rho^n {\bf v}\right > $
then exists the third and formula \ref{main_f1}  is valid.
This can be simply shown using lemmas 2.2,2.3.

\vskip 2em

%\pagebreak

{\bf Corollary 3.1}

Let ${\bf v}, {\bf w}$ - the elements of degree 1, from $S V$ and $S V^*$ respectively,
i.e. ${\bf v}\in V, {\bf w}\in V^*$   then:

\bea\frac{Tr_{S V}\left (\bar\rho A (e^{\bf w}) C (e^{\bf v})\right)
}{Tr_{S V}\left (\bar\rho\right)}
=e^{\left < {\bf w}|\frac{1}{1 -\rho} {\bf v}\right >}
\label{main_sl_exp}\eea

One can prove corollary expanding exponents to the series and applying
theorem 3.1, or one can obtain straightforward proof of corollary 3.1,
using commutation relations $ e^{A({\bf w})}e^{C({\bf v})}=$ $e^{<{\bf w}|{\bf v}>} e^{C({\bf v})} e^{A({\bf w})}$,
and reasoning in the same way, as it was done in proof 1 of the theorem 2.1.

\vskip 1em

Proofs of theorem 3.1 are similar to the ones of theorem 2.1.
One should only substitute sign plus for sign minus in proofs 1 and 3 of
theorem 2.1 to obtain the ones of theorem 3.1.
One can also change proof 2 of theorem 2.1 in order to obtain proof of
theorem 3.1. But a little difference arises in this case. Let us consider it.

$SV= $ $\CC[e_1,e_2,....e_N]=\CC[e_1]\otimes\CC[e_2]\otimes...\otimes\CC[e_N]$,
where $e_i$ any basis of $V$. Hence, similar as it was done in previous
section one can reduce the proof to the case $dim V=1$.
But $dim SV= \infty$, in contrast to $dim \L V= 2$. So the proofs
differs for the case $dim V=1$.
Let us prove theorem 3.1 for the case $dim V=1$.

Let $x \neq 0 \in V$; ~~ $ \rho x = \l x$; ~~
hence: $\bar\rho (x^n) = \l^n x^n$.
Let $\partial_x \in V^*$ be dual to $x$.
\bea
&&Tr_{SV}(\bar \rho A(\tilde x^k) C(x^k) )= Tr_{\CC[x]} \bar \rho \frac{\partial^k}{\partial x^k } x^k =
\sum_{n=0}^{\infty} \frac{(n+k)!}{k!} (\l)^n=\frac{\p^k}{\p {\l}^{k}} (\sum_{n=0}^{\infty} (\l)^n)=
\nn\\
&&
=\frac{\p^k}{\p {\l}^{k}} \frac{1}{1-\l}
=\frac{k!}{(1-\l)^{k+1}}=
(Tr_{S V} \bar\rho) k!<\tilde e | \frac{1}{1-\rho} e>^k=
(Tr_{S V} \bar\rho) <(\tilde e)^k | \frac{1}{1-\bar\rho} e^k>\nn
\eea
We have used obvious fact that $Tr_{SV} \bar\rho = \frac{1}{1-\l}$.
Note that one can rewrite it in the form: $det_{V}\frac{1}{(1-\rho)}= Tr_{SV} \bar\rho$.

It is also evident that
if $k\neq l$, then
$Tr_{\CC[e]} \bar\rho (\partial_e)^k e^l =0=<(\tilde e)^k | \frac{1}{1-\bar\rho} e^l>$.

So the case $dim V=1 $ is completed. Hence, proof 2 is finished.

\vskip 0.5em

We have seen that $det_{V}\frac{1}{(1-\rho)}= Tr_{SV} \bar\rho$ in case $dim V=1$.
Similar as it was done in previous section one can prove it for any $V$.
Hence the following proposition is true:

{\bf Proposition 3.2:} Let $\rho$ be an operator with trace on
space $V$ and $|\rho|<1 $, then
$$det_V \frac{1}{1-\rho}=\sum_{n=0}^{\infty}Tr_{S^n V} (\bar \rho )$$

\vskip 0.5em

Let us formulate analogs of propositions 2.3-6

{\bf Proposition 3.3}: Let $e_i$ be arbitrary basis of
space $V$, $\tilde e_i$ be the dual basis,
 then for any operator O on space $S^n V$ the following formula
 takes place:

\bea Tr_{S^n V} O=\frac{1}{n!}\sum_{l_1... l_n\geq 0} {<\tilde e_{l_1}
\tilde e_{l_2}...\tilde e_{l_n}|O e_{l_1} e_{l_2}... e_{l_n} > }&&
\label{sim_sim}\eea

{\bf Corollary 3.2:} Let  $\rho$ be an operator on finite-dimensional space $V$,
$\tilde e_i$ be the dual basis, $\bar\rho$ is induced operator on $S^n V$.
$e_i$ - basis of $V$. Then the following formula takes place:
\bea
Tr_{S^n V} \bar \rho= \frac{1}{n!}
\sum_{l_1...l_n \geq 0}
per \left\vert\matrix{
<\tilde e_{l_1} |\rho e_{l_1} >& <\tilde e_{l_2} |\rho e_{l_1} >&...&<\tilde e_{l_n} |\rho e_{l_1} > \cr
<\tilde e_{l_1} |\rho e_{l_2} >& <\tilde e_{l_2} |\rho e_{l_2} >&...&<\tilde e_{l_n} |\rho e_{l_2} > \cr
 \vdots &  \vdots & \ddots & \vdots  \cr
<\tilde e_{l_1} |\rho e_{l_n} >& <\tilde e_{l_2} |\rho e_{l_n} >&...&<\tilde e_{l_n} |\rho e_{l_n} >
}\right\vert
&&
\label{sim_sim_det_2}
\eea

{\bf Proposition 3.4} Let $\rho$ be an operator with trace on space $V$, $
\bar\rho$ is induced operator on $S^n V$, then $Tr_{S^n V}\bar\rho $
can be expressed through $Tr_{V}\rho^l $ as follows:

\bea Tr_{S^n V}\bar\rho=\sum_{ n_j\geq0:\sum_j jn_j=n}
\frac{1}{n_1! n_2!... 1^{n_1}2^{n_2}...}\prod_{l} (Tr_V\rho^l) ^{n_l}
\eea

{\bf Proposition 3.5} Let $\rho $ be an operator with trace on space $ V$; $ {\bf v_i}, {\bf w_i}$
arbitrary elements from $V$ and $V^*$ respectively. Then trace of
 $\bar\rho\prod_{i=1}^k {A ({\bf w_i})}\prod_{j=1}^k {C ({\bf v_j}) }$ over
space $S^nV$ can be  expressed through
$Tr_{S^r V}\rho$ and $ < {\bf w_i}|~\bar\rho ^{l}~{\bf v_j} > $ as follows:

\bea Tr_{S^n V}\left (\bar\rho\prod_{i=1}^k {A ({\bf w_i})}
\prod_{j=1}^k {C ({\bf v_j})}\right)=
\sum_{\begin{array}{c} r, q_1, q_2,..., q_k\geq0:\\r+\sum_j q_j=n
\end{array}} Tr_{S^r V} (\bar\rho)\left <\prod_{j=1}^k {\bf w_j}~|~
\prod_{j=1}^k\rho^{q_j} ({\bf v_j})\right > &&\label{lem}
\eea

The proofs of the propositions above are similar to the ones in previous
section.

Due to corollary 3.2 one can reformulate proposition 3.2 as follows:

{\bf Proposition 3.6:} Let $\rho$ be an operator with trace on
finite-dimensional space $V$, $e_i$ a basis of $V$,
$\tilde e_i$ be the dual basis.

\bea
det_{V} \frac{1}{(1-\rho)}=Tr_{S V} \bar\rho= \sum_{n=0}^{\infty}
\frac{1}{n!} \sum_{l_1,l_2,...,l_n}
per | <\tilde e_{l_j} |\rho e_{l_i}>|_{i,j\leq n}
\label{per_2}
\eea

%%%%%%%%%%%%%%%%%%%%%%%%%%%%%%%%%%%%%%%%%%%%%%%%%%%%%%%%%%%%%%%%%%%%%%%%%%
%%%%%%%%%%%%%%%%%%%%%%%%%%%%%%%%%%%%%%%%%%%%%%%%%%%%%%%%%%%%%%%%%%%%%%%%%%
%%%%%%%%%%%%%%%%%%%%%%%%%%%%%%%%%%%%%%%%%%%%%%%%%%%%%%%%%%%%%%%%%%%%%%%%%%
%%%%%%%%%%%%%%%%%%%%%%%%%%%%%%%%%%%%%%%%%%%%%%%%%%%%%%%%%%%%%%%%%%%%%%%%%%

\setcounter{equation}{0}

\section{Fredholm's formulas }

In this section we recall definitions of Fredholm's  determinant
and minor, and also Fredholm's   formula  for the solution
of integral equation. We show that Fredholm's determinant
coincides with determinant defined  above (definition 2.3).
We also show that Fredholm's minor   can be expressed,
as the trace of certain operator over the space $\L V$. We show
that our main formula (theorem 2.1) in the case,
when $\rho$ is integral operator, $deg {\bf v}=deg{\bf w}=1$ coincides
with Fredholm's formula   for the solution of integral equations. Thus
we obtain simple  proof of Fredholm's formula, since
proof 1 of the theorem 2.1 occupies only few lines.
Also we shall show that the trace formula for the bosonic Fock space
(theorem 3.1)
leads to the formula of inversion of integral operator, but
in terms of the "Fredholm's permanents".

Let $V=C [a, b] $ ( precisely speaking one must consider
$C^* [a,b] $; look remark 4.1).
Consider integral equation $\phi (x)+\int K (x, y)\phi (y)=f (x) $.
Let $\rho$ be an integral operator with kernel $K (x, y) $.

Recall the definitions of  Fredholm's determinant   and Fredholm's minor:
\bea
D^{fred} (1+\rho) &=& 1+\int K (\xi,\xi) d\xi+\frac{1}{2!}\int
\int d\xi_1d\xi_2 det\left\vert\matrix{ K (\xi_1,\xi_1) & K (\xi_1,
\xi_2)\cr
                           K (\xi_2,\xi_1) & K (\xi_2,\xi_2)}\right
\vert\nn\\&&+...+\frac{1}{n!}\int\int...\int d\xi_1d\xi_2... d\xi_n
det| K (\xi_i,\xi_j) |_{i, j\geq n}+.......\\D_{s, t}^{fred} (1+\rho)
& =& K (s, t)+\int d\xi_1 det\left\vert\matrix{ K (s, t) & K (s,\xi_1)
\cr
                           K (\xi_1, t) & K (\xi_1,\xi_1)}\right\vert
\nn\\&&+\frac{1}{2!}\int\int d\xi_1d\xi_2 det\left\vert\matrix{ K (s,
t) & K (s,\xi_1) & K (s,\xi_2)\cr
                           K (\xi_1, t) & K (\xi_1,\xi_1) & K (\xi_1,
\xi_2)\cr K (\xi_2, t) & K (\xi_2,\xi_1) & K
                           (\xi_2,\xi_2)}\right\vert
+.......
\eea

Considering the set $\delta(x-s)$, at $s\in [a,b]$ as  "continuous basis" in
$C [a, b] $ one can obviously see that   the definition of the Fredholm's
determinant  is absolutely analogous to the  formula \ref{det_2}.
Hence it is naively clear that $D^{fred}(1+\rho)$ $=Det(1+\rho)$
= $Tr_{\L V} \bar\rho$.
Let us prove it rigorously.

The action of an integral operator with continuous kernel can be canonically
extended to the action on delta-functions:
$\int K(x,y) \delta(x-s) dx = K(s,y)$, obtained function is conituous.

{\bf  Lemma 4.1} Let $\rho$ be an integral operator with
continuous kernel $K(x,y)$, then
\bea
Tr_{\L^n V} \bar\rho=
\frac{1}{n!} \int\int...\int d\xi_1d\xi_2... d\xi_n
<\delta(x-\xi_1)\land...\land \delta(x-\xi_n)|\bar\rho
\delta(x-\xi_1)\land...\land \delta(x-\xi_n)>
\eea
Proof: consider the case  $n=1$  hence we need to prove that:
$Tr \rho= \int d\xi K(\xi, \xi)$. This fact is naively obvious.
In order to prove it rigorous we note that any continuous
function $K(x,y)$ can be represented in the form:
$K(x,y)=\sum_i \l_i f_i(x)g_i(y)$, where $\sum |\l_i| <\infty$,
$|f_i(x)|<C, |g_i(y)|<C$, $f_i,g_i$ are continuous functions.
By definition 2.2 trace $\rho$ equals to $\sum_i \l_i \int d\xi
f_i(\xi)g_i(\xi)$, hence it's clear that $Tr \rho=\int d\xi K(\xi, \xi)$.
In order to prove the formula for the case $n>1$ one should use
the formula \ref{yvno}. Lemma is proved.

Using lemma one can easily obtain that
$D^{fred}(1+\rho)=$ $Det(1+\rho)$ = $Tr_{\L V} \bar\rho$.

Analogically Fredholm's minor    equals to:

\bea D_{s, t}^{fred} (1+\rho) =Tr_{\L V} A (\delta (x-s))\overline{
\rho} C (\delta (x-t))
\label{minor1}
\eea

Actually:
\bea &&Tr_{\L^n V} A (\delta (x-s))\bar\rho C (\delta (x-t))
=\frac{1}{n!}\int...\int d\xi_1.... d\xi_n <\delta (x -\xi_1)\land...
\land\delta (x -\xi_n)||A (\delta (x-s))\nn\\&&\bar\rho C (\delta (x-t))
\delta (x -\xi_1)\land...\land\delta (x -\xi_n) >
=\frac{1}{n!}\int...\int d\xi_1.... d\xi_n <\delta (x-s)\land\delta (x -
\xi_1)\land...\land\delta (x -\xi_n)|\nn\\&&|\rho (\delta (x-t))
\land\rho (\delta (x -\xi_1))\land...\land\rho (\delta (x -\xi_n)) >
\nn\\&&
=\frac{1}{n!}\int...\int d\xi_1.... d\xi_n det\left\vert\matrix{ K (s,
t) & K (s,\xi_1) &... & K (s,\xi_n)\cr K (\xi_1, t) & K (\xi_1,\xi_1)
                           &... & K (\xi_1,\xi_n)\cr
\vdots &\vdots &\ddots &\vdots &\cr K (\xi_n, t) & K (\xi_n,\xi_1) &....
                           & K (\xi_n,\xi_n)}\right\vert
\nn\eea

In order to make the reasonings above rigorous, let us note three
facts. The first, $\delta(x-s)\notin C[a,b]$, hence there is problem
in defining $C(\delta(x-s))$, but due to $K(x,y) $  is continuous
function the product $\rho C(\delta(x-s))$ is well defined.
The second,  one can naturally extend the action of $ C(\delta(x-s))$ and
$A(\delta(x-t))$ on the delta-functions: $A(\delta(x-t))[\delta(x-r)]=
\delta_{t,r}$, $C(\delta(x-t))[\delta(x-r)]=\delta(x-t)\land \delta(x-r)$.
The third,  operator
$A( \delta(x-s) ) \overline{ \rho} C(\delta(x-t))$ on the space
$\L C[a,b]$ can be represented
in the form: $\sum_i \rho_{i<L} + \sum_{i<K} T_i$, where $\rho_i$ are operators
with continuous kernels, $T_i$ are operators of the form: $T_i= t_i(x)\otimes
 \delta(x-t_i)$. Where  $t_i$ are continuous functions, $\delta(x-t_i)$
are considered as functional on $C[a,b]$. Hence that operator is operator
with trace. And it's easy to see that lemma 4.1 holds for this operator also,
if action of this operator is defined on delta-functions as we have just
described.

%\newpage

One can prove lemma 4.1 in another way.

Really, expression

$
\frac{1}{n!} \int\int...\int d\xi_1d\xi_2... d\xi_n
<\delta(x-\xi_1)\land...\land \delta(x-\xi_n)|\bar\rho
\delta(x-\xi_1)\land...\land \delta(x-\xi_n)>
$
can be transformed, in the same way, as we have done
proving the proposition  2.4. But one
should write integrals instead of sums. After the transformations one obtains:

$\sum_{ n_j\geq0: \sum_j jn_j=n }
\frac{1}{n_1!n_2!...1^{n_1}2^{n_2}...}
\prod_{l} (-1)^{(l-1)n_l}(\int...\int d\xi_1...d\xi_l K(\xi_1,\xi_2)K(\xi_2,\xi_3)....K(\xi_l,\xi_1))^{n_l}$.

Applying Merser's theorem ($Tr_{V} \rho=\int d\xi K(\xi,\xi)$).
One gets:

$\sum_{ n_j\geq0:\sum_j jn_j=n}
\frac{ (-1) ^n}{n_1! n_2!... 1^{n_1}2^{n_2}...}\prod_{l} (-Tr_V (\rho)
 ^l) ^{n_l}$ that is equal to $Tr_{\L^n V}\bar\rho$
according to proposition 2.4. Proof is finished.

One can similar prove the formula \ref{minor1} using
the proposition 2.5 instead of 2.4.

Thus we obtained the following proposition:

{\bf Proposition 4.1:}

\bea D^{fred} (1+\rho) &=&Tr_{\L V}\bar\rho=Det (1+\rho)\\
D_{s,t}^{fred} (1+\rho) &=&Tr_{\L V} A (\delta (x-s))\overline{\rho} C (
\delta (x-t))\eea

Hence, theorem 3.1 gives us expression for the
matrix  elements of the operator $\frac{1}{1+\rho}$:

\bea <\delta (x-s)|\frac{\rho}{1+\rho}\delta (x-t) >
=\frac{Tr_{\L V} A (\delta (x-s))\overline{\rho} C (\delta (x-t))}
      {Tr_{\L V}\bar\rho}=
\frac{D_{s, t}^{fred} (1+\rho) }{D^{fred} (1+\rho)} &&\nn\\<\delta (
x-s)|\frac{1}{1+\rho}\delta (x-t) >=<\delta (x-s)|\delta (x-t) >
\frac{D_{s, t}^{fred} (1+\rho) }{D^{fred} (1+\rho)}\label{rho_inv}&&\eea

Thus, one obtains Fredholm's formula   for the solution
of integral equation $\phi (x)+\int K (x, y)\phi (y)=f (x) $:

\bea\phi (s) =f (s) -\int dt\frac{D_{s, t}^{fred} (1+\rho) }{D^{fred}
(1+\rho)} f (t)\label{solut}\eea
Really, $\phi (x)=\frac{1}{1+\rho}f (x) $,
hence, due to formula \ref{rho_inv} $\phi (s) =f (s) -
\int dt\frac{D_{s, t}^{fred} (1+\rho) }{D^{fred} (1+\rho)} f (t) $

Remark 4.1: Simpler to obtain formula \ref{solut}
using not the $Tr_{\L V}\bar\rho A (\delta (x-s))C (\delta (x-t)) $, but
the trace: $Tr_{\L V}\bar\rho A(\delta (x-s)) C (f (x)) $.
Handling this way, we shall obtain the same formula,
but we would not leave the space $C[a,b]$, as it happened when we
consider $Tr_{\L V}\bar\rho A (\delta (x-s))C (\delta (x-t)) $
(note that $\delta(x-t) \notin C[a,b]$).
But we wanted to emphasize analogy between our formulas  and Fredholm's ones,
so we have admitted ourselves some inaccurateness.

\vskip 2em

Theorem 3.1 is analog of theorem 2.1 for bosonic Fock space.
Using it one can obtain formulas similar to \ref{solut}.
But expressions involves not determinants, but permanents.

Let us consider
an integral equation $\phi (x) - \int K (x, y)\phi (y)=f (x) $,
at condition $|\rho|^n\to 0 $. (Where $\rho$ is an integral operator
with kernel $K (x, y)$.)

Let us define $P^{fred}(1-\rho)$ and $P_{s,t}^{fred} (1 -\rho)$ as follows:
\bea P^{fred} (1 -\rho) &=&\sum_{n=0}^{\infty}\frac{1}{n!}\int\int...
\int d\xi_1d\xi_2... d\xi_n ~~per|K (\xi_i,\xi_j) |_{i, j\leq n}=
\sum_{n=0}^{\infty} Tr_{S^n V}\bar\rho=det_{V}\frac{1}{1 -\rho}\\
P_{s,t}^{fred} (1 -\rho) & =&\sum_{n=0}^{\infty}\frac{1}{n!}\int\int...
\int d\xi_1d\xi_2... d\xi_n ~~per\left\vert\matrix{ K (s, t) & K (s,
\xi_1) &... & K (s,\xi_n)\cr
                           K (\xi_1, t) & K (\xi_1,\xi_1) &... & K (
\xi_1,\xi_n)\cr\vdots &\vdots &\ddots &\vdots &\cr
                           K (\xi_n, t) & K (\xi_n,\xi_1) &.... & K (
\xi_n,\xi_n)}\right\vert\nn\\&=&
\sum_{n=0}^{\infty} Tr_{S^n V} A (\delta (x-s))\bar\rho C (\delta (x-t)
\eea

Hence theorem 3.1 gives us the formulas  for matrix  elements of
the inverse operator:

\bea <\delta (x-s)|\frac{\rho}{1 -\rho}\delta (x-t) >
=\frac{Tr_{S V} A (\delta (x-s))\overline{\rho} C (\delta (x-t))}
      {Tr_{S V}\bar\rho}=
\frac{P_{s, t}^{fred} (1 -\rho) }{P^{fred} (1 -\rho)} &&\nn\\<\delta
(x-s)|\frac{1}{1 -\rho}\delta (x-t) >=<\delta (x-s)|\delta (x-t) >+
\frac{P_{s, t}^{fred} (1 -\rho) }{P^{fred} (1 -\rho)}
\label{rho_inv_2}&&\nn\eea

Hence, we obtain the analog of Fredholm's formula   for the solution
of integral equation $\phi (x) -\int K(x, y)\phi (y)=f (x) $:

\bea\phi (s) =f (s)+\int dt\frac{P_{s, t}^{fred} (1 -\rho) }{P^{fred}
(1 -\rho)} f (t)\eea

Remark: maybe, the  name "Fredholm's permanents"
is not acquitted, because $P^{fred} (1 -\rho) =Tr_{S V}\bar\rho =det_{V}
\frac{1}{1 -\rho}$. The last equality is Proposition 3.6.

%%%%%%%%%%%%%%%%%%%%%%%%%%%%%%%%%%%%%%%%%%%%%%%%%%%%%%%%%%%%%%%%%%%%%%%%%%
%%%%%%%%%%%%%%%%%%%%%%%%%%%%%%%%%%%%%%%%%%%%%%%%%%%%%%%%%%%%%%%%%%%%%%%%%%
%%%%%%%%%%%%%%%%%%%%%%%%%%%%%%%%%%%%%%%%%%%%%%%%%%%%%%%%%%%%%%%%%%%%%%%%%%
%%%%%%%%%%%%%%%%%%%%%%%%%%%%%%%%%%%%%%%%%%%%%%%%%%%%%%%%%%%%%%%%%%%%%%%%%%

\setcounter{equation}{0}

\section{Traces of intertwining operators.}

In this section  we apply  general formulas obtained above for the
the computing the traces of some
concrete operators, which encounters during the solution
of integrable models.
Proposition 4 is used in paper \cite{KLP2} for the computing
formfactors of $SU(2)$-invariant Thirring model. It easily follows from
the Proposition 2. Formula from the  proposition 2 was suggested by
the author's scientific advisor S.M. Khoroshkin. And the  proof
of this formula was the initial purpose of this paper.
Further the author came to the general formulas (theorems 2.1,3.1).

Actaully, traces which we calculate
in this section, turn out to be divergent that
is why it is necessary to regularize them somehow.
We carried out the discussion of
regularization in appendix 2. Also
on this moment it's not known explicit
topological description of the Fock space where the operators acts.
But it's clear that this space is subspace of the
$\CC [[a_1, a_2,...]] $ and monomials $ (a_1) ^{i_1}
(a_2) ^{i_2}... (a_k) ^{i_k}$ are the basis of this space, in some sense
of the word "basis".
That is why one may understand these operators
as infinite matrixes, and the  trace as the sum of diagonal elements,
our theorem is valid and in this situation, this
is discussed in appendix 1. In this section, we do not discuss
all these problems but formally apply theorem 3.1 and its corollary 3.1.

Let $V$ be a linear space.

Let $\aa (u),\b (v) $ be  functions with values in $V^*, V $ respectively,
such that pairing between them is equal $g (u-v) $.

Let operator $e^{\g d}$ acts on $V$ so that its pointwise action on
functions can  be written as follows:
$e^{\g d}\b (z)=\b (z+\g) ~~~~~~ (e^{\g d}) ^*\aa (z)=\aa (z -\g) $.

{\bf Proposition 5.1.}

\bea\frac{1}{Tr_{S V} e^{\g d}} Tr_{S V} e^{\g d} e^{\aa (u)} e^{\b (
v) }=\prod_{n=0}^{\infty} e^{ g (u-v-n\g)}\label{gg}\eea

As it have been already said, the   traces  of similar type are used
for the solution of integrable models. This proposition is useful,
because it provides the expression for the trace in terms of pairing only.

One can obtain the proof by simply applying corollary 3.1.

\vskip 1em

Let $ a_{n}\in V, a_{-n}\in V^* $ such elements that $ < a_{-m}| a_{n} > =n\delta_{n, m}$.

Let $\tilde a_+ (u)=\sum_{i=1}^{\infty}\frac{a_{i}}{i} u^{-i} ~~~~
\tilde a_ - (v)=\sum_{i=1}^{\infty}\frac{a_{-i}}{i} v^{i}$.

Obviously $ <\tilde a_ - (u),\tilde a_+ (v) >=-log (1 -\frac{u}{v}) $.

{\bf Proposition 5.2.}

\bea\frac{Tr_{SV} (e^{\gamma d}\prod_{i=1}^{I}\exp (k_i\tilde a_ - (
\alpha_i))\prod_{j=1}^{J}\exp (l_j\tilde a_+ (\beta_j)) }{Tr_{SV} (e^{
\gamma d})}=\prod_{m=1}^\infty\prod_{i, j} (\beta_j -\alpha_i - m
\gamma) ^{-k_i l_j}\label{formul2}\eea where $\sum k_i =0\sum l_j=0 $.

Proof easily follows from corollary 3.1.

Note that obtained infinite product is convergent. Really, this
follows from following proposition:
\beq\prod_{k_1,\ldots, k_n\geq0} {\prod_m (a_m+\sum k_j\omega_j)\over
\prod_p (b_p+\sum k_j\omega_j)}\label{Barnes}\eeq Is alighted when:
\beq\sum_m (a_m) ^q=\sum_p (b_p) ^q,\qquad q=0, 1,\ldots, n
\label{constr}\eeq (We prove this proposition in appendix 5).

And the conditions of this proposition are fulfilled, since
$\sum k_i =0\sum l_j=0 $.

\vskip 1em

Actually, in physical applications one uses not the space of polinomials
$SV=\CC [a_1, a_2,....] $, but the sum  of the infinite numbers of the copies
of this space: $F=\bigoplus_{M\in\ZZ/2} F^M$, where $F^M={\CC} [[a_{-1},
\ldots, a_{-n},\ldots]]\otimes\left ({\CC}e^{Ma_0}\right) $.

Let us define the action  of operator $p$ as follows:

\bea p (a_{-1}^{i_1}a_{-2}^{i_2}.... a_{-m}^{i_m}\otimes e^{na_0})=2n
  (a_{-1}^{i_1}a_{-2}^{i_2}.... a_{-m}^{i_m}\otimes e^{na_0})
\label{p}\eea

Let us define the following functions:
$ a_+ (u)=\sum_{i=1}^{\infty}\frac{a_{i}}{i} u^{-i}-plog
(u) ~~~ a_ - (v)=\sum_{i=1}^{\infty}\frac{a_{-i}}{i} v^{i}+
\frac{a_0}{2} $

Let an operator $e^{\g d}$ acts on the space  $F$ so that:

\bea e^{\g d} a_ - (v) =a_ - (v+\g)\label{gamma_d}\eea

{\bf Proposition  5.3:} Let us define operator $O$ on space $F$ as follows:\\
$O=\prod_{j=1}^{J}\exp (l_j a_+ (\beta_j))\prod_{i=1}^{I}\exp (k_i a_ -
(\alpha_i))\mbox{ where }\sum k_i =0\sum l_j=0 $

Then takes place the following  formula for his   trace:
\bea\frac{Tr_{F} e^{\g d} O }{Tr_{F} (e^{\gamma d})}=\prod_{m=0}^
\infty\prod_{i, j} (\beta_j -\alpha_i - m\gamma) ^{-k_i l_j}
\label{formul333}\eea

{\bf Proof: }

It is obvious that due to condition $\sum k_i =0$,
each subspace $F^M\subset F$ is invariant under the action of operator $O$.

Due to proposition 2 :
\bea\frac{Tr_{F^0} e^{\g d} O }{Tr_{F^0} (e^{\gamma d})}=\prod_{m=1}^
\infty\prod_{i, j} (\beta_j -\alpha_i - m\gamma) ^{-k_i l_j}\nn\eea

Let us show that operator $O e^{\g d}$ on the subspace $F^M$ is
equivalent to the operator $\prod_j\b_j^{2M l_j} O$ $e^{2M\tilde a_ - (\g)} e^{\g d}$
on subspace $F^0$,i.e.  corresponding matrix elements are equal.
\bea
&& < e^{Ma_0}\otimes w~ O e^{\g d} v\otimes e^{M a_0} >=<
e^{Ma_0}\otimes w~ O [(e^{\g d} v)\otimes e^{\g d} e^{M a_0}] >=\nn\\
&&
\mbox{From formula \ref{gamma_d} follows that:} e^{\g d}
a_0=a_0+2\tilde a_-(\g)\mbox{ hence:}\nn\\
&&
=< e^{Ma_0}\otimes w~ O [e^{2M\tilde a_ - (\g)} (e^{\g d} v)]\otimes O
[e^{M a_0}] >=\nn\\&&
\mbox{from  formula \ref{p} follows that:}
 O [e^{M a_0}]=\prod_j\b_j^{2M l_j} e^{M a_0}ю\mbox{ hence:}\nn\\
&&
=\prod_j\b_j^{2M l_j} < e^{Ma_0}\otimes w~ O [e^{2M\tilde a_ - (\g)} (
e^{\g d} v)]\otimes e^{M a_0}] >
=\prod_j\b_j^{2M l_j} < 1\otimes w~ O e^{2M\tilde a_ - (\g)} e^{\g d} v
\otimes 1] >\nn
\eea
Thus, we have obtained, in particular, the following proposition:

Lemma: $ Tr_{F^M} e^{\g d} O=\prod_j\b_j^{2M l_j} Tr_{F^0} e^{\g d} O
e^{\tilde a_ - (\g)} $

By the direct computation, using Corollary 3.1, and
making use of the conditions:$\sum k_i=\sum l_j=0$ we obviously get
$$\prod_j\b_j^{2M l_j} Tr_{F^0} e^{\g d} O e^{
\tilde a_ - (\g)}=(Tr_{F^0} e^{\g d})\prod_{m=1}^\infty\prod_{i, j} (
\beta_j -\alpha_i - m\gamma) ^{-k_i l_j} $$ Hence, expression
in LHS is not dependent on M.

Thus, we see that   traces  of operators $e^{\g d} O, e^{\g d}$ over
all spaces $F^M$ are the same. Thus, their quotient,
just like it was required in proposition 3 equals to:
\bea\prod_{m=1}^\infty\prod_{i, j} (\beta_j -\alpha_i - m\gamma) ^{-
k_i l_j}\label{rrrrrr6}\eea
Therefore,  the quotient of the  traces of these operators over the space
$F=\bigoplus_M F^M$ is equal to:\ref{rrrrrr6}.

{\bf Proof is finished. }

\vskip 1em

Let us  introduce the following generating function:
\beq\eta_+ (z)=\lim{K\to\infty} (2\h K) ^{ -\da/2}
\prod_{k=0}^K e^{a_+ (z-2k\h) - a_+ (z -\h-2k\h)}\label{eta}\eeq
It can easily be checked that it's well-defined.

In paper \cite{KLP1} it was shown that
intertwining operators for basic representations of central extended
Yangian double of $\frak{sl}_2$ can be written as follows:
$ e^{\epsilon_1 a_ - (\aa)}\eta_+^{
\epsilon_2} (\aa\pm\h) $, where $\epsilon_i=\pm 1$

In paper \cite{KLP2} the formfactors  of $SU(2)$-invariant Thirring model were expressed
through the traces of the products of intertwining operators.
In that paper the following formula have been used:

{\bf Proposition 5.4:}

\bea\frac{\tr_{F}\left (e^{\gamma d}\prod_j e^{k_j a_ - (w_j)}\prod_k
\eta_+^{p_k} (z_k)\right) }{tr_F e^{\g d}}=\prod_{m=1}^\infty
\prod_{k=0}^\infty\prod_{j, k} {(z_k-w_j -\h-m\g-2k\h) ^{k_jp_k}\over (
z_k-w_j-m\g-2k\h) ^{k_jp_k}}\nn \label{calc-gener}
~~~~\mbox{ where }\sum_j k_j=0\sum_k p_k =0\label{cond-gener}\nn
\eea
Due to \ref{Barnes} and $\sum_j k_j=0 ~ \sum_k p_k =0 $
the product in RHS of the equality is convergent.

To prove the proposition one should substitute
the definition of $\eta_+ (z)$ (\ref{eta}) to the proposition 5.3.

\setcounter{equation}{0}
\section{Acknowledgements.}

The author thanks his  scientific adviser S. M. Khoroshkin for the
formulation of the task, significant help during the work  and
constant support. Also author thanks S. Pakulyak for multiple explanations
and B. Feigin who pointed to the \cite{SJM}.
The research of the author was supported by RFBR grant 96-15-96939
and INTAS grant 93-0183.

%%%%%%%%%%%%%%%%%%%%%%%%%%%%%%%%%%%%%%%%%%%%%%%%%%%%%%%%%%%%%%%%%%%%%%%%%%
%%%%%%%%%%%%%%%%%%%%%%%%%%%%%%%%%%%%%%%%%%%%%%%%%%%%%%%%%%%%%%%%%%%%%%%%%%
%%%%%%%%%%%%%%%%%%%%%%%%%%%%%%%%%%%%%%%%%%%%%%%%%%%%%%%%%%%%%%%%%%%%%%%%%%
%%%%%%%%%%%%%%%%%%%%%%%%%%%%%%%%%%%%%%%%%%%%%%%%%%%%%%%%%%%%%%%%%%%%%%%%%%

%\newpage

\app{}

\sapp{~Space with basis.}

In main text we have proved theorems for the normed spaces.
In this appendix we shall show that if one consider
infinite-size matrixes then theorems are also valid, at some conditions
on matrixes.

Note that one can consider infinite-size matrixes as operators
acting: $\CC^\infty_0~~\rightarrow~~ \CC^\infty$. Where
$\CC^\infty_0$ is space of finite sequences,
$\CC^\infty$ is the space of all sequences.
Certainly,  the product of two infinitesize matrixes is not always well-defined.

The trace of infinitesize matrix  is sum of its diagonal elements.

Let $ {\bf v},{\bf w} \in  \CC^\infty $. Define pairing as follows:
$<{\bf v}|{\bf w}>=\sum_i v_iw_i $

We shall say that trace, product, applying matrix on vector, pairing are
WELL-DEFINED, iff the series arising in this operations are absolutely
convergent.

One defines algebras $S \CC^\infty_0$, $S \CC^\infty$,
$\L \CC^\infty_0$, $\L \CC^\infty$, in the same way as it was done
above for normed $V$.
For any ${\bf v} \in \CC^\infty$ one defines creating-annihilating operators
$A({\bf v}), C({\bf v})$ acting:$\CC^\infty_0~~\rightarrow~~ \CC^\infty$.

Theorems 2.1 and  3.1 are true in this case:

{\bf Theorem}
Let $\rho$ be an infinitesize matrix, ${\bf w}, {\bf v} \in \L \CC^\infty$.
Assume that $\rho^n, Tr  \rho^n , <{\bf w}|\rho^n {\bf v}>  $,
$<{\bf w}|\frac{1}{1+\rho} {\bf v}>$ are well-defined.
Then
$Tr_{\L \CC^\infty_0 \rightarrow \L \CC^\infty} \bar \rho A({\bf w}) C({\bf v})$
and  $Tr_{\L \CC^\infty_0 \rightarrow \L \CC^\infty} \left( \bar \rho \right)$,
are well-defined too and the following formula takes place:
\bea
\frac{Tr_{\L \CC^\infty_0 \rightarrow \L \CC^\infty} \left(
\bar \rho A({\bf w}) C({\bf v})  \right)}
{Tr_{\L \CC^\infty_0 \rightarrow \L \CC^\infty} \left( \bar \rho \right)}
=\left< {\bf w}|\frac{1}{1+\rho} {\bf v} \right>
\label{in_f2_2}
\eea

{\bf Theorem}
Let $\rho$ - be an infinitesize matrix, ${\bf w}, {\bf v} \in S \CC^\infty$.
Assume that $\rho^n, Tr  \rho^n , <{\bf w}|\rho^n {\bf v}>  $
are well-defined and series
$\sum_i Tr \rho^i ,~~~\sum_n <{\bf w}|\rho^n {\bf v}>$
are absolutely convergent.
Then
$Tr_{S\CC^\infty_0 \rightarrow S \CC^\infty} \bar \rho A({\bf w}) C({\bf v})$
and  $Tr_{S\CC^\infty_0 \rightarrow S \CC^\infty} \left( \bar \rho \right)$
are well-defined and the following formula takes place:
\bea
\frac{Tr_{S\CC^\infty_0 \rightarrow S \CC^\infty} \left(
\bar \rho A({\bf w}) C({\bf v})  \right)}
{Tr_{S\CC^\infty_0 \rightarrow S \CC^\infty} \left( \bar \rho \right)}
=\sum_{n=0}^{\infty} \left< {\bf w}| \bar\rho^n {\bf v} \right>
\label{in_f2_3}
\eea

To prove these theorems it suffices to note that proof 3 is valid
in this situation, because the only difference is that here arises
series instead of finite sums, but the series are absolutely convergent.

\vskip 2em

\sapp{ Regularization of the traces.}

It is possible that in theorem 3.1 LHS of equality  does  not exist
while RHS exists. That is why naturally to give the following definition.

{\bf Definition:}

\bea Reg\frac{Tr_{SV}\left (\bar\rho A ({\bf w}) C ({\bf v})\right)
}{Tr_{SV}\left (\bar\rho\right)}
=\left\{\matrix{\left < {\bf w}|\frac{1}{1 -\rho} {\bf v}\right > ~~
\mbox{  }\nn\\
\sum_{n=0}^{\infty}\left < {\bf w}|\rho^n {\bf v}
\right >}\right.\nn
\eea
We choose that value in RHS which exists, if both of them exists then
they are obviously equal to each other.

\vskip 1em
According to proposition 3.2 $ Tr_{SV}\bar\rho =det_{V}\frac{1}{1 -\rho}$
when  $|\rho| < 1$. That is why naturally to introduce the following
definition.

{\bf Definition:} $Reg Tr_{SV}\bar\rho=det_{V}
\frac{1}{1 -\rho}$ if $\frac{1}{1 -\rho}$ exists.

The  traces of such type are used for computing formfactors, in
quantum field theory and statistical physics. In  quantum field
theory models traces usually turns out to be divergent,
though  we do not have rigorous proof that our regularization
coincides with scaling limit from corresponding lattice model
or satisfies Smirnov's axioms (\cite{S1}) on formfactors or agreed with
ultraviolet cut-off (\cite{L}), but the concrete example
of computing for $SU(2)$-invariant Thirring model shows that
all approaches lead to the  same results.

Further we shall show that in some cases such regularization coincides
with certain limit of the quotient of   traces, and in another
cases   with regularization by means of zeta function (if the one exists,
though, for $SU(2)$-invariant Thirring model or Sine-Gordon model such regularization through
zeta function does not exist.)

{\bf regularization by passage to the limit.}

Let us consider the
operator $\rho=e^{\g d}$ that  was used in section 5.
In natural basis $a_1,a_2,....$ it has the form:
$$\left (\matrix{ 1 & 0 & 0 &...\cr
\g & 1 & 0 &...\cr\frac{\g^2}{2!} &\g & 1 &...\cr
\vdots &\vdots &\vdots &\ddots}\right) $$
and obviously  $Tr_{V}\rho=\infty$, hence
$Tr_{SV}\bar\rho A ({\bf w}) C ({\bf v})=\infty, ~~~ Tr_{SV}\bar\rho=\infty$.
But nevertheless, for concrete ${\bf w}, {\bf v}$ that were used
in section 5, quotient of these  traces was equal to $\sum_{i} < {\bf w}|
\rho^i ({\bf v}) $ and was convergent. Further on we shall
give yet another argument in favour that, in this
case  definition $Reg\frac{Tr_{SV}
\left (\bar\rho A ({\bf w}) C ({\bf v})\right) }{Tr_{SV}\left (\bar
\rho\right)}
=\sum_{n=0}^{\infty}\left < {\bf w}|\rho^n {\bf v}\right > $ is acquitted.

Let us introduce the operator $T_t$ such that
$ T_t(a_i)=(t)^ia_i$.
Obviously, at $0<t<1$ operators  $T_t$, $\overline{T_t}$ are operators with
trace on the spaces $V$ and $SV$, respectively.
$Tr_{V} T_t= \frac{1}{1-t} $, ~~~ $Tr_{V} T_t= \prod_{i}\frac{1}{1-t^i} $.
Also, it is clear that $T_t\rho$ is operator with trace and
$(T_t \rho)^i \to 0 ~~as  i\to \infty$, hence,
$\overline{T_t\rho}$ is operator with trace on $SV$.
Diagonal elements of  $\overline{T_t \rho}$ in basis
$a_1^{k_1}a_2^{k_2}a_3^{k_3}...$
tends to diagonal elements of $\bar\rho$ as $t\to1$, hence,
it's natural to consider
$Reg_{(lim)}~~~ \frac{Tr_{SV} \bar\rho A({\bf w}) C({\bf v})}{Tr_{SV} \bar \rho}=
\lim{t\to 1}
\frac {Tr_{SV} \left( T_t\bar \rho A({\bf w}) C({\bf v})  \right)}
{Tr_{SV}\left( T_t \bar \rho \right)}$

{\bf Proposition:} $\lim{t\to 1}\frac
{Tr_{SV}\left (T_t\bar\rho C ({\bf w}) C ({\bf v})\right)} {Tr_{SV}
\left (T_t\bar\rho\right) }=\sum_{n=0}^{\infty}\left < {\bf w}|\bar
\rho^n {\bf v}\right >\label{st_reg2}$ i.e. $Reg_{ (lim) }~~~
\frac{Tr_{SV}\bar\rho A ({\bf w}) C ({\bf v}) }{Tr_{SV}\bar\rho}= Reg
\frac{Tr_{SV}\bar\rho A ({\bf w}) C ({\bf v}) }{Tr_{SV}\bar\rho}$

The proof is based on Abel's theorem:
if $\sum_{i} a_i <\infty $, then $\lim{t\to1}\sum_i a_i t^i=\sum_i a_i$.
We will not go into details.

Thus we have obtained that regularization through the limit coincides
with our main regularization.

Let us generalize the example described above.

Let $V$ be a normed vector space with the Grothendieck's approximation property.
It's known that there exists the sequence of operators with trace $H_n$
that converge uniformly on any compact subspace to identical operator.
It's clear that one can  choose $H_n$ such  that
$\lim{i\to\infty} (H_n) ^i=0$.

{\bf Proposition:} Assume that  $\lim{i\to
\infty} (\rho H_n) ^i =0$ and series  $\sum_{n=0}^{
\infty}\bar\rho^n {\bf v}$ is convergent. Then $Tr_{SV} H_n\bar\rho <
\infty$, $Tr_{SV} H_n\bar\rho A ({\bf w}) C ({\bf v}) <\infty$ and

\bea\stackreb{\mbox{\rm lim}}{n\to\infty}\frac {Tr_{SV}\left (H_n\bar
\rho C ({\bf w}) C ({\bf v})\right)} {Tr_{SV}\left (H_n\bar\rho\right)
}=\sum_{n=0}^{\infty}\left < {\bf w}|\bar\rho^n {\bf v}\right >
\label{st_reg1}\eea

{\bf Proof.}

Proof is based on the observation that one can apply the theorem to
 to operators $H_n\rho$, because $\lim{i\to\infty} (\rho H_n) ^i =0$.
Hence:
\bea\stackreb{\mbox{\rm lim}}{n\to\infty}\frac {Tr_{SV}\left (H_n\bar
\rho C ({\bf w}) C ({\bf v})\right)} {Tr_{SV}\left (H_n\bar\rho\right)}
=\stackreb{\mbox{\rm lim}}{n\to\infty}\sum_{n=0}^{\infty}\left < {\bf
w}| (H_n\bar\rho) ^n {\bf v}\right >=\sum_{n=0}^{\infty}\left < {\bf
w}|\bar\rho^n {\bf v}\right >
\eea
Last passage to the limit is acquitted, because  $\bar\rho^n {\bf v} $
is convergent, hence the set    $\bar\rho^n {\bf v} $ ~~$n>0\in\ZZ$
is precompact set, consequently $H_n\to\id$ is convergent uniformly
on this set, so on can pass to the limit.

{\bf Proposition is proved.}

Note that, apparently, such regularization is not the limit
of quotient  of traces over finite-dimensional subspaces exhausting
every bit of $V$. If $Tr_{SV}\bar\rho A (
{\bf w}) C ({\bf v})=\infty ~~~Tr_{SV}\bar\rho=\infty$ then,
apparently, the limit of their quotient over finite-dimensional subspaces is
is equal $\infty$ too.

For example consider the case: $dim V=1$,  $\rho=id$,

$\lim{N\to\infty}\frac{\sum_{i=0}^N Tr_{S^i V}\bar\rho A ({\bf w}) C (
{\bf v}) }{\sum_{i=0}^N Tr_{S^i V}\bar\rho}=\infty$

\vskip 2em

{\bf Regularization by means of  the analytical continuation.}

Regularization by means of passage to the limit
was possible, when $\rho$ look like on identical operator.
In the case if $|\rho| > 1$ then regularization through the passage to
the limit is impossible.

It's possible another variant of regularization by means of analytical
continuation, in particular by means of zeta function, we shall show,
that this variant coincides with with our main regularization.
Though note that in the previous example such method of regularization
is not applicable.

Consider $\zeta (s) =Tr_{SV}\bar\rho^{ (-s) }$, assume that  for
$s$ large enough $Tr_{SV}\bar\rho^{ (-s) }$ is defined and assume that  this
function can be  continued to meramorphic on the entire plane.
Then we shall define $Reg_{\xi}~~Tr_{SV}\bar\rho=\zeta (1) $.
Note that typically $Tr_{V}\rho^s $ has  pole in unit. But
$Tr_{SV}\bar\rho^{-s}=Det_V\frac{1}{1 -\rho^s} $
 and this function is analityc in point $s=1$, if $1\notin Spec~~\rho$.
Similar: we can define  $\tilde\zeta (s) =
Tr_{SV}\bar\rho^{-s} A ({\bf w}) C ({\bf v}) $,
and put  $Reg_{\xi} Tr_{SV}\bar\rho A ({\bf w}) C ({\bf v})=\tilde\zeta (1) $.

Let $s$ be sufficiently large, then:

\bea Reg_{\xi}\frac{Tr_{SV}\left (\bar\rho^s A ({\bf w}) C ({\bf v})
\right)} {Tr_{SV}\left (\bar\rho^s\right) }=\frac{Tr_{SV}\left (\bar
\rho^s A ({\bf w}) C ({\bf v})\right)} {Tr_{SV}\left (\bar\rho^s\right)
}=\left < {\bf w}|\frac{1}{1 -\rho^s} {\bf v}\right >
\eea
Function $\zeta (s) $ is analitycal, hence,  if equality is true
at  $s > N$ that it is true and always, i.e.
$Reg_{\xi}\frac{Tr_{SV}\left (\bar\rho A ({\bf w}) C ({\bf v})\right)}
{Tr_{SV}\left (\bar\rho\right) }=\frac{\tilde\zeta (1) }{\zeta (1) }=
\left < {\bf w}|\frac{1}{1 -\rho}
{\bf v}\right > =Reg\frac{Tr_{SV}\left (\bar\rho A ({\bf w}) C (
{\bf v})\right)} {Tr_{SV}\left (\bar\rho\right) }$
Thus  regularization by means of zeta function coincides with our.

The author thinks that the same situation will be
with any other analitycal regularization, because
we always want to deform our operator somehow to the area, where traces
are defined, continue the function from this area to the hole plain.
But in the area, where traces are defined our formula is true, hence
it will be true on the hole plane, because if two analitycal functions are
equal in some area , they are equal on the hole plain.

\sapp{ Continuos basis.}

 Considering the models of the quantum field theory,
for example Sine-Gordon model \cite{KLP3}, \cite{JKM},
it's necessary to consider
% $SV=C [a_1, a_2,...] $ to consider
free bosons $a_t$, where $ t\in\RR$. $a_t$  plays
role of "continuous basis" in space $V$.
Instead of giving rigorous definition of continuous basis, we prefer
to consider two examples. The first is the following:
any function belongs  $L_2 (\RR) $ one can expand to the Fourier's integral:
$f (x)=\int\tilde f (\l) e^{i\l x}d\l$ i.e. $e^{i\l x}$ at $\l\in\RR$ is a
continuous basis $L_2 (\RR) $; the second: for any continuous
function $f (x)$ holds that $f (x)=\int f (y)\delta (y-x) dy $ i.e.
$\delta (y-x) $ at $x\in [a, b] $ is  continuous basis in $C [a, b] $.
Common in this two examples is the following:  basic functions do not belong
the space, but the elements of the space
can be expressed as follows: $\int d\l h (\l) e_{\l}$, where $e_{\l} $ is
continuous basis, $h(\l) $ is certain class of scalar functions
determining the  space.

 In the case of continuously indexed bosons there must exist similar
description of space $V$. Space $V$ is the set of the following elements:
$\int d\l h (\l) a_{\l}$, where $h (\l) $ belongs certain class of scalar functions.
Usually $V$ is a Hilbert space, respectively, our theorem will be valid
in this case.

If $a_t$ is continuos basis of $V$, then trace of any operator $\rho$
on $V$ can be found by the formula:
$Tr_V \rho = \int\limits_{\l >0} <a_{-\l} |\rho a_{\l}>d\l$,
in case when integral is convergent.

Respectively, the trace over space $SV$,  can be found
by the formula similar to proposition 3.3:

$Tr_{SV} O = \sum_{n=0}^{\infty}
\int\limits_{\l_1\leq 0}... \int\limits_{\l_n\leq 0} d\l_1...d\l_n
\frac{ <a_{-\l_1}a_{-\l_2}....a_{-\l_n}
|O a_{\l_1}a_{\l_2}....a_{\l_n}>}{ <a_{-\l_1}a_{-\l_2}....a_{-\l_n}
|a_{\l_1}a_{\l_2}....a_{\l_n}>}$

One can consider this equalities as definitions of traces.
Theorem 2.1, 3.1 are valid, because one can apply proof 3, substituting
sums for integrals.

\sapp{ Traces and generating functions.}

In this appendix we shall introduce the "by-product" of our work:
initially the proof of our main results (theorems 2.1, 3.1)
used the propositions below, but  then we understood that they are
not necessary. But may be it will be interesting for somebody.

As it is known, considering affine algebras it's convenient to work
with generating functions i.e. with functions with values in algebra or
the space of its representation. For example, it's easier to
write down the  commutation relations in terms of  generating functions,
but not the concrete generators.

 Initially, our purpose was to prove the  concrete proposition (\ref{formul2}).
In this proposition we considered ${\bf v}=a_{-} (\aa) ,{\bf w}=a_{+} (\b)$.
$a_{\pm} (z) $ are the generating functions, the action of
operator $\rho = e^{\g d} $ can be easily written in terms of
generating functions, but rather complicated in terms of concrete generators.
Therefore it was natural to try to express the trace in terms of generating functions.

 So the main purpose of this section to find expressions for the traces
in terms of generating functions. Turns out to be,
that it's possible to  each operator on space $V$ assign some integral operator
acting on generating functions.
Generating functions plays role of the "continuous basis"
and the trace can be expressed as some integral.

 We shall only formulate the propositions, we shall not prove them,
because the proofs are trivial.

Let $V$ be a  linear space, $dim V=N\leq\infty$. $V^{*}$ -
its dual. Let $v_i~~ 0\leq i\leq N-1 $ be a  basis
in $V$; $w_j~~ 1\leq j\leq N $ - basis in $V^*$.

Let us denote by: $\b (x)=\sum_{j=1}^{N} w_j (x) ^{-j}$~~ ~~ $
\aa (y)=\sum_{i=0}^{N-1} v_i (y) ^{i}$.

Let us introduce "identical pairing tensor":

\bea d (x, y)=\frac{1}{x}\sum_{i=0}^{N-1}\frac{y^i}{x^i}\eea

{\bf Proposition  A. 4.1:} Vectors $v_i, ~~~w_j $ are the dual
basis in $V, V^{*}$ respectively, iff  $ <\b (x)|\aa(y) > =d (x, y) $.

{\bf Definition:  } We shall define the formal integral from power series $F(x)$
as follows: $ (\int F (x) dx) $ equals to  coefficient at $\frac{1}{x}$.

{\bf Proposition  A. 4.2} Let $A$ be an operator on space $V$.
The function: $\widetilde{ A} (x, y)=\sum_{0\leq i, j\leq N-1}
a^i_j\frac{z^j}{x^{i+1}}$ corresponds to $A$,
 where $a^i_j$ are matrix  elements of the
operator $A$ in basis $v_i$. Then the following  formulas holds:

\bea\widetilde{ (\id)} (x, y) &=&d (x, y)\\A (\aa (z)) &=&\int
\widetilde A (x, z)\aa (x) dx\\\label{pr2_2}\widetilde{AB} (x, y) &=&
\int\widetilde B (x, z)\widetilde A (z, y) dz\\Tr A &=&\int A (x, x) dx\\
\widetilde{\sum_i w_i\ot v_i} (x, y) &=& <
\b (x)|\aa (y) >\label{pr2_5}\eea

Remark 1:
Formal integral defined above, coincides with
ordinary integral, if we choose the appropriate contour of integration.
In the case $dim V <\infty$ one can
choose arbitrary contour which includes the zero and not includes $\infty$.
In case $dim V=\infty$ the functions  $\aa (z),\b (z),\widetilde{A} (x, y) $
may have some poles, therefore one must choose contour with
respect to them.
General principle is the following: typically we consider
integrals from the following expressions:
for the integral: $\int F_1 (z) F_2 (z) dx $,  where $F_1 (z) $ is the
formal power series of $\frac{1}{z}$; ~~~ $F_2 (z) $ is
formal power series of $z$, contour must contain singular points of
$F_1(z)$ and does not contain singular points of $F_2(z)$.
In examples under consideration only integrals of such type arise.
For example in formula \ref{pr2_2} contour must be selected as follows:
at each  $x$ contour must contain all singular points for $A (x, z) $
and does  not contain all singular points for $\aa (z) $.
In proposition A.4.3 contour must contain
singular points for $\b (z) $ and does  not contain singular points for
$A (\aa (z)) $. Note that
that singular points for functions $ A(\aa (z)) $ and $\aa (z) $
may be different, since, operator $ A $ is not obliged to be
bounded. Also we want to note that broadly speaking, it's   not necessarily to require
the analytical continuation  of functions $ A\left (\aa (z)\right),\b (z) $
to the hole plain, it's sufficient to require that areas where they not analitycal one can separate
by  a contour, and integrate over it.

Remark 2: In this section we mean that $Tr A$ is sum of diagonal elements in
basis $v_i$. In the case $dim V=\infty$ we shall consider only such
operators $A$ that sum of their diagonal elements is absolutely convergent.

{\bf Proposition  A. 4.3} Let us assume that $ <\b (x)|\aa (y) > =d (x, y) $,
then for the trace of any operator $A$ over space $V$ takes place the formula:

\bea Tr_{V} A=\int <\b (x) |A\aa (x) > dx\label{dual_fields_fml}\eea

Remark: if $e_i$, $e^i$ are  dual basis
in $V$ and $V^{*}$ respectively, then for any operator $A$ holds that:
$Tr A=\sum_i < e^i A e_i > $. It's clear that
formula \ref{dual_fields_fml} is the analogous to  above formula of linear
algebra, if one  considers $\aa (z),\b (z) $ as dual basis
indexed by $ z $. One may be surprised
why we call by  dual basis, such basis  that
 pairing between them  is equal $d (u, v)=\frac{
\frac{u^{N}}{v^{N}}-1}{v-u}$, but not the delta-function $\delta (u, v) $.
But $d (u, v) $ may be considered as
$\delta (u, v) $, since
 $\int d (u, v) P (u) ~du=P (v) $ for any polinomial $P(u) $
of degree not exceeding $N-1$. This property is similar to the
one of delta-function. One can see that this is the only
property, which is necessary.

\vskip 1em

{\bf Proposition  A.4.4.} Let us denote by $g (x, y)=<\b (x)|\aa (y) > $.
Due to the property \ref{pr2_5} $g(x,y)$  is the "kernel" of operator
$\sum_i w_i\ot v_i$. Let $g^{-1} (x, y) $ be  the "kernel" of the inverse
operator. Note,  that integral operators with "kernels" $g (x, y) $ and $g^{-1} (x, y)$
 mutually inverse.

\bea Tr_V A=\int\int g^{-1} (x, y) <\b (x) |A\aa (x) > dx dy\eea

Remark: this Proposition is
the analogue  of following simple  fact of linear algebra:
$Tr A=\sum_{i, j} (g^{-1}) ^{i}_{j} < w^{j} A v_i > $, where
$g^{i}_{j}= < w^{i}|v_j > $.

\vskip 1em

In some cases $g^{-1} (x, y) $ can be found
by the following simple method.

{\bf Proposition A.4.5} Let us that $g (x, y)=A_y d (x, y) $. Where ${A_y }$
is  self-adjoint operator acting on functions from variable
$y$. Then $g^{-1} (x, y)={A_y}^{-1} d(x, y) $

Proof:

\bea\int g^{-1} (z, y) g (x, y) dy & =&\int {A_y}^{-1} d (z, y) {A_y}
 d (x, y) dy=\nn\\
=\int d (z, y) {A_y}^{-1} {A_y} d (x, y) dy &=&\int d (z, y) d (x, y)
dy=d (z, y)\nn\eea

%From[Out of] call 2 easily to see that
%proven property[trait] is necessary[is it is necessary] for that, in
%order $g^{-1} (x, y)={A_y}^{-1} d (x, y) $. But in effect[as a matter
%of fact|in actual fact|as a fact|as it is] it and it's sufficient,
%because if $\int K (u, v) g (u, y)=0 $ and function K (u, v) exist
%polinomiayl degrees[extents|grades|rates|scales] not exceeding $N$ from
%$u,\frac{1}{v}$ that K (u, v) =0. Integral operators with kernels
%[nuclei] g (x, y) and $g^{-1} (x, y) $ mutually inverse.

Similar one can obtain:

{\bf Proposition A.4.6}: Let us assume that
$g (x, y)={A_x} d (x, y) $. Where ${ A_x}$ is self-adjoint
operator acting on functions
from variable $x$. Then $g^{-1} (x, y)={A_x}^{-1} d (x, y) $

\vskip 1em

Similar as it have been done above, one
can obtain the following propositions,
analogous to lemmas 2.1, 3.1.

{\bf Proposition  A.4.7} Let us assume that
 $ <\b (x)|\aa (y) > =d (x, y) $, then for the
trace of any operator $A$ over the space $S^nV$, $\L V$ respectively,
takes place the formulas:

\bea Tr_{S^nV} A=\frac{1}{n!}\int...\int <\b (x_1)\b (x_2)...\b (x_n)
|A\aa (x_1)\aa (x_2)...\aa (x_n) >\prod_{i=1}^{n} dx_i\label{pr7} &&\\
Tr_{\Lambda^n V} A=\frac{1}{n!}\int...\int <\b (x_1)\land\b (x_2)
\land...\land\b (x_n) |A\aa (x_1)\land\aa (x_2)\land...\land\aa (x_n) >
\prod_{i=1}^{n} dx_i
 &&\eea

{\bf Proposition A. 4.8}
Let us assume that $ <\b (x)|\aa (y) > =g (x, y) $, then for the
trace of any operator $A$ over the space  $S^nV$, $\L V$
respectively, takes place the formulas:

\bea Tr_{S^nV} A=\frac{1}{n!}\int...\int\prod_{i=1}^{n} g^{-1} (x_i,
y_i) <\prod_{i=1}^{n}\b (x_i) |~A~\prod_{i=1}^{n}\aa (y_i) >
\prod_{i=1}^{n} dx_idy_i\label{pr8} &&\\Tr_{\Lambda^n V} A=\frac{1}{n!}
\int...\int\prod_{i=1}^{n} g^{-1} (x_i, y_i) <\bigwedge_{i=1}^n\b (x_i)
|~A~\bigwedge_{i=1}^n\aa (y_i) >\prod_{i=1}^{n} dx_idy_i &&\eea

\sapp{Convergence of infinite product.}

{\bf Proposition A.5.1 }
If
\beq\sum_m (a_m) ^q=\sum_p (b_p)^q,\qquad q=0, 1,..., n
\label{constr2}\eeq
then the following infinite-product is convergent:
\begin{equation}\prod_{k_1,..., k_n\geq0} {\prod_{m=1}^{M} (a_m+
\sum k_j\omega_j)\over\prod_{p=1}^{P} (b_p+\sum k_j\omega_j)}\\
\label{Barnes2} ~~~~~~Where ~\Re\omega_i < 0\nn\\
\eeq

Apparently, Barns \cite{Ba} was the first who studied
similar products in his theory of generalized gamma-functions.
The proof desribed below was communicated to the author by his
scientific adviser S.M. Khoroshkin.

{\bf Proof:}

Note,
that M=P (otherwise n-th the member of product will not tends to 1).

Proof is based on integral representation for \ref{Barnes2}.\\
Let us recall Frullani's formula:
if integral $\int_A^{\infty}\frac{f (x) }{x}$ exists $\forall A > 0$, then
$\int_0^{\infty}\frac{f (ax)-f (bx) }{x}=\ln (b/a) f(0) $.

Hence:
\beq
\ln({\prod_{m=1}^{N} (a_m+ k\omega)\over\prod_{p=1}^{N}
(b_p+k\omega) })=\int_0^{\infty}\frac {\sum_{m=1}^{N} e^{ (a_m+k
\omega) x} -\sum_{p=1}^{N} e^{ (b_p+k\omega) x}}{x}
\eeq

We can consider $a_i, b_i < 0$. Since
for the large enough number $k_j$ holds: $ (a_m+\sum
k_j\omega_j) ~ and  ~ (b_p+\sum k_j\omega_j)  < 0$.

Let us prove at first that
\begin{equation}\prod_{k\geq0} {\prod_{m=1}^{M} (a_m+k\omega)\over
\prod_{p=1}^{M} (b_p+\sum k\omega) }=\nn\\\exp (\int_{0}^{\infty}
\frac {\sum_{m=1}^{M} e^{a_m x} -\sum_{p=1}^{M} e^{b_p x}}{ (x) (1 -
\exp (\omega x)) })\eeq where $\Re\omega < 0$

Really:

\bea
\ln (\prod_{k\geq0} {\prod_{m=1}^{N} (a_m+ k\omega)\over
\prod_{p=1}^{N} (b_p+k\omega) })=\sum_{k\geq0}\ln ({\prod_{m=1}^{N} (
a_m+ k\omega)\over\prod_{p=1}^{N} (b_p+k\omega) })=\nn\\\sum_{k
\geq0} (\int_{0}^{\infty}\frac {\sum_{m=1}^{N} e^{ (a_m+k\omega) x} -
\sum_{p=1}^{N} e^{ (b_p+k\omega) x}}{x}=\nn\\\int_{0}^{\infty}
\frac {\sum_{m=1}^{N}\left ((e^{a_m x}-e^{b_m x})\sum_{k\geq0} e^{k
\omega x}\right) }{x}=
\int_{0}^{\infty}\frac {\sum_{m=1}^{N} (e^{a_m x}-
e^{b_m x})}{(1- \exp (\omega x))x}\nn
\eea

Similar

\bea\prod_{k_1,\ldots, k_n\geq0} {\prod_{m=1}^{M} (a_m+\sum k_j
\omega_j)\over\prod_{p=1}{P} (b_p+\sum k_j\omega_j) }=\exp\int_{0}^{
\infty}\frac {\sum_{m=1}^{N} (e^{a_m x}-e^{b_m x}) }{x\prod_{i=1}^{n} (
1 -\exp (\omega_i x))}\eea

Hence, the convergence of the product is equivalent
to the convergence of the integral. The integral is
obviously convergent in the point   $\infty$.
Conditions of convergence in point zero give us
\begin{equation}
\nn\sum_m (a_m) ^q=\sum_p (b_p) ^q,\qquad q=0, 1,\ldots, n\eeq

%
%%%%%%%%%%%%%%%%%%%%%%%%%%%%%%%%%%%%%%%%%%%%%%%%%%%%%%%%%%%%%%%%%%%%%%%%%%%
% Список литературы
%%%%%%%%%%%%%%%%%%%%%%%%%%%%%%%%%%%%%%%%%%%%%%%%%%%%%%%%%%%%%%%%%%%%%%%%%%%%
%\newpage
%\newcommand{\bib}[4]{\bibitem{#1}#2, ``{\em#3\/}'', #4}

\end{document}